\documentclass{article}

\usepackage{arxiv}

\usepackage[utf8]{inputenc} 
\usepackage[T1]{fontenc}    
\usepackage{hyperref}       
\usepackage{url}            
\usepackage{booktabs}       
\usepackage{amsfonts}       
\usepackage{nicefrac}       
\usepackage{microtype}      
\usepackage{lipsum}
\usepackage{graphicx}
\graphicspath{ {./images/} }

\usepackage{color,soul}

\usepackage{amsmath,amssymb}
\DeclareMathOperator*{\argmax}{arg\,max}
\DeclareMathOperator*{\argmin}{arg\,min}
\usepackage{multirow}
\newcommand{\tcr}{\textcolor{red}}

\title{Optimal Decision Making in\\High-Throughput Virtual Screening Pipelines}

\author{
  Hyun-Myung~Woo\\
  Computational Science Initiative\\
  Brookhaven National Laboratory\\
  Upton, NY 11973, USA\\
  \texttt{hwoo@bnl.gov}\\
  \And
  Xiaoning~Qian\\
  Department of Electrical and Computer Engineering\\
  Texas A\&M University\\
  College Station, TX 77843, USA\\
  \texttt{xqian@ece.tamu.edu}\\
  \And
  Li~Tan\\
  Computational Science Initiative\\
  Brookhaven National Laboratory\\
  Upton, NY 11973, USA\\
  \texttt{ltan@bnl.gov}\\
  \And 
  Shantenu~Jha\\
  Computational Science Initiative\\
  Brookhaven National Laboratory\\
  Upton, NY 11973, USA\\
  \texttt{shantenu@bnl.gov}\\
  \And 
  Francis~J.~Alexander\\
  Computational Science Initiative\\
  Brookhaven National Laboratory\\
  Upton, NY 11973, USA\\
  \texttt{falexander@bnl.gov}\\
  \And 
  Edward~R.~Dougherty\\
  Department of Electrical and Computer Engineering\\
  Texas A\&M University\\
  College Station, TX 77843, USA\\
  \texttt{edward@ece.tamu.edu}\\
  \And
  Byung-Jun~Yoon\\
  Department of Electrical and Computer Engineering\\
  Texas A\&M University\\
  College Station, TX 77843, USA\\
  \texttt{bjyoon@ece.tamu.edu}\\
}

\begin{document}
\maketitle
\begin{abstract}
The need for efficient computational screening of molecular candidates that possess desired properties frequently arises in various scientific and engineering problems, including drug discovery and materials design. However, the large size of the search space containing the candidates and the substantial computational cost of high-fidelity property prediction models makes screening practically challenging. In this work, we propose a general framework for constructing and optimizing a virtual screening (HTVS) pipeline that consists of multi-fidelity models. The central idea is to optimally allocate the computational resources to models with varying costs and accuracy to optimize the return-on-computational-investment (ROCI). Based on both simulated as well as real data, we demonstrate that the proposed optimal HTVS framework can significantly accelerate screening virtually without any degradation in terms of accuracy. Furthermore, it enables an adaptive operational strategy for HTVS, where one can trade accuracy for efficiency.

\end{abstract}



\section{Introduction}
\label{sec:introduction}
In various real-world scientific and engineering applications, the need for screening a large set of molecular candidates to identify a small subset of molecules that satisfy certain criteria or possess targeted properties arises fairly frequently. For example, since the Coronavirus disease $2019$ (COVID-$19$) outbreak, there have been significant concurrent efforts among various groups of scientists to identify or develop drugs that can provide a potential cure for this extremely infectious disease. One such notable effort is IMPECCABLE (Integrated Modeling PipelinE for COVID Cure by Assessing Better LEads)~\cite{saadi2020impeccable} whose operational objective is to optimize the number of promising ligands that potentially lead to the successful discovery of drug molecules. To this aim, IMPECCABLE utilized deep learning-based surrogates for predicting docking scores and multi-scale biophysics-based computational models for computing docking poses of compounds. Built on the strength of massive parallelism on exascale computing platforms combined with RADICAL-Cybertools (RCT) managing heterogeneous workflows, IMPECCABLE identified promising leads targeted at COVID-$19$.

Considering that severe acute respiratory syndrome coronavirus $2$ (SARS-CoV-$2$), the virus that can potentially lead to COVID-$19$, is known to rapidly mutate itself to create more infectious and deadlier variants~\cite{roy2021global}, such drug screening process to identify anti-viral drug candidates that are effective against a specific variant may have to be repeated as new variants emerge. However, when one considers the huge search space of potential molecules--\textit{e.g.}, ZINC (Zinc Is Not Commercial) $15$~\cite{sterling2015zinc} contains about $230$ million commercially available compounds. There are, however, about ${10^{12}}$ compounds that can be considered for drug design in chemical space theoretically--and the astronomical amount of computation that was devoted to the screening of drug candidates in~\cite{saadi2020impeccable} to screen ${10^{11}}$ candidates, this is without question a Herculean task that requires enormous resources and one that cannot be routinely repeated.

While different in scale and complexity, such high-throughput virtual screening (HTVS) pipelines have been widely utilized in various fields, including biology~\cite{rieber2009rnaither,studer2010engineering,hartmann2011htpheno,sikorski2018high,clyde2021high}, chemistry~\cite{martin2014silico,cheng2015accelerating,chen2017developing,filer2017tcpl,rebbeck2020ryr1,saadi2020impeccable}, engineering~\cite{tran2020smf}, and materials science~\cite{yan2017solar,zhang2020first}.
However, the construction of such HTVS pipelines and the strategies for operating them heavily rely on expert intuition, often resulting in heuristic methods with reasonable yet sub-optimal screening performance. 
It remains a fundamental challenge to optimally construct and operate such screening pipelines to sift potential molecular candidates from an enormous search space in an efficient yet accurate manner.

\begin{figure*}[h!]
\centerline{\includegraphics[width=\textwidth]{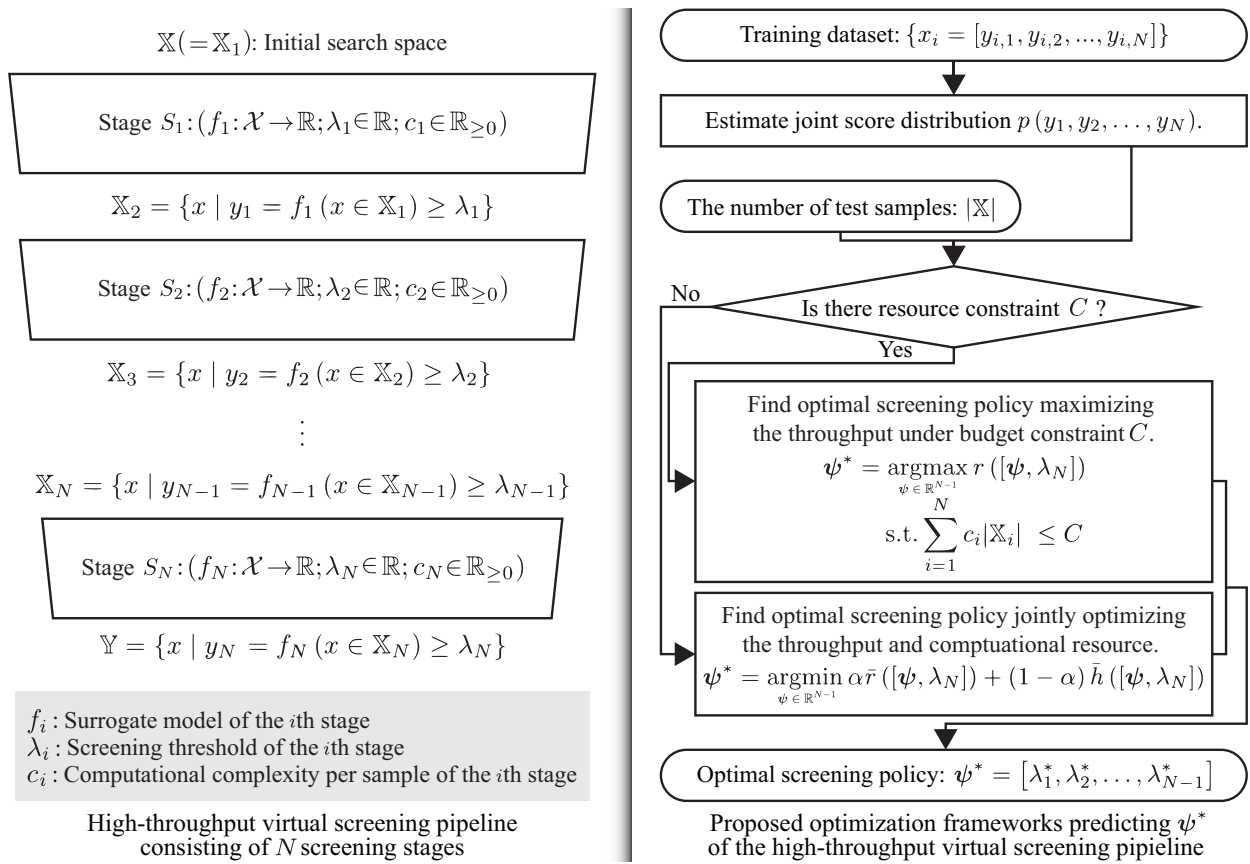}}
\caption{Illustration of a general high-throughput virtual screening (HTVS) pipeline (left) that consists of $N$ stages (surrogate models) for rapid and reliable identification of a set ${\mathbb{Y}}$ of candidate molecules that likely possess the desired properties from a huge original set ${\mathbb{X}}$ containing all candidates. Stage ${S_i}$ evaluates all the molecules ${x \in \mathbb{X}_i}$, which passed the previous stage ${S_{i-1}}$, via a surrogate model ${f_i}$. ${S_i}$ passes the sample ${x}$ to the next stage ${S_{i+1}}$ if ${f_i\left(x\right) \ge \lambda_i}$. Otherwise, it discards the molecule. The proposed optimization framework shown on the right side predicts optimal screening policy ${ \boldsymbol{\psi}^\ast  = \left[\lambda_1^\ast, \lambda_2^\ast, \dots, \lambda_{N-1}^\ast \right]}$ that yields the maximal throughput according to the screening campaign scenarios.}
\label{fig1}
\end{figure*}

In general, typical HTVS pipelines consist of multiple stages, each of which is associated with a surrogate model that evaluates the property of the molecules with a different accuracy/fidelity and computational cost. This is illustrated on the left side of Fig.~\ref{fig1}. At each stage in the pipeline, the molecular candidate is evaluated to determine whether the evaluation result appears promising enough to warrant passing it to the next--often more computationally expensive but more accurate--stage without unnecessarily wasting computational resources and time. In this way, the HTVS pipeline narrows down the number of candidate molecules, while sensibly allocating the available resources for investigating those that are promising and more likely to possess the desired property. The most promising candidates that remain at the end of screening may proceed to experimental validation, which is often more laborious, costly, and time-consuming. For example, in~\cite{clyde2021high}, an HTVS pipeline based on multi-fidelity surrogate models combined with an experimental platform successfully selected and reported a novel non-covalent inhibitor, MCULE-5948770040. The reported inhibitor has been identified by screening over $6.5$ million molecules, and it has been shown to inhibit the SARS-Cov-$2$ main protease. HTVS pipelines have been also widely used for materials screening. For example, a first-principles high-throughput screening pipeline for non-linear optical materials (FHSP-NLO)~\cite{zhang2020first} consisting of several computational predictors, based on density functional theory (DFT) calculations as well as data transformation and extraction methods, successfully identified deep-ultraviolet non-linear optical crystals that were reported in previous studies~\cite{chen1995new,shi2017finding,zhang2017fluorooxoborates,luo2018m2b10o14f6,mutailipu2018srb5o7f3,wang2018cation,zhang2018cab5o7f3}. 

Although previous studies have demonstrated the advantages of constructing an HTVS pipeline for rapid screening of huge set of molecules to narrow down the most promising molecular candidates that are likely to possess the desired properties, the problem of \textit{optimal decision-making} in such screening pipelines has not been extensively investigated to date. For example, how should one decide whether or not to pass a molecular candidate at hand to the next stage, given the outcome of the current stage? More specifically, in the HTVS example shown on the left side of Fig.~\ref{fig1}, how do we optimally determine the screening threshold of each stage for a given HTVS structure? Furthermore, if we were to modify the HTVS structure or construct it from scratch by interconnecting multi-fidelity surrogate models, what would be the optimal structure of such HTVS that maximizes the throughput and accuracy? This requires selecting the optimal subset of the available multi-fidelity models, arranging them in the optimal order, and then exploring the interrelations among their predictive outcomes to make optimal operational decisions for the constructed HTVS.

In this paper, we present a computational framework that can answer the aforementioned questions and applied to the optimization of HTVS pipelines involving that consist of multiple surrogate models with different costs and fidelity. The key idea is to estimate the joint probability distribution of predictive scores that result from the different stages constituting the HTVS pipeline, based on which we optimize the screening threshold values. We consider two optimization scenarios. First, we consider the case where the total computational budget is fixed and the goal is to maximize the throughput within the given budget. Second, we consider the case where we aim to jointly maximize the throughput of the HTVS pipeline while minimizing the overall computational cost required for screening. We demonstrate the performance of the proposed HTVS pipeline optimization framework based on both \textit{simulated data} as well as \textit{real data}. In the simulated example, the joint distribution of the predictive scores from the multi-fidelity models at different stages is assumed to be known, based on which we extensively evaluate the performance of the proposed approach under various scenarios. As a second example, we consider the problem of screening for long non-coding RNAs (lncRNAs). In this example, we first construct an HTVS pipeline by interconnecting existing lncRNA prediction algorithms with varying costs and accuracy and apply our proposed framework for performance optimization. Both examples clearly demonstrate the advantages of our proposed scheme, which leads to a substantial reduction of the total computational cost at virtually no degradation in overall prediction accuracy. Furthermore, we show that the proposed framework enables one to make an informed decision to balance the trade-off between speed and accuracy, where one could trade accuracy for higher efficiency, and vice versa.

\section{Methods}
\label{sec:methods}
In this section, we first formally describe the operational process of a general HTVS pipeline and provide a high-level overview of the proposed optimization framework. Two different scenarios will be considered. In the first scenario, the objective is to maximize the expected throughput of the HTVS pipeline under a fixed computational budget. In the second scenario, the objective is to jointly optimize the screening throughput and the computational efficiency of the pipeline.


\subsection*{Overview of the proposed HTVS pipeline optimization framework}
We assume that an HTVS pipeline consists of ${N}$ screening stages ${S_i: \left(f_i: \mathcal{X} \rightarrow \mathbb{R}; \lambda_i; c_i \right)}$, ${i=1, 2, \dots, N}$, connected in series as shown in Fig.~\ref{fig1} (left), where ${f_{i}: \mathcal{X} \rightarrow \mathbb{R}}$ is a surrogate model for predicting the property of interest for a given molecule and ${\lambda_i}$ is the screening threshold. The average computational cost per sample for ${f_i}$ associated with the ${i}$th stage ${S_i}$ is denoted by ${c_i}$.
At each stage ${S_i}$, the corresponding surrogate model $f_{i}$ is used to evaluate the property of all molecules ${x \in \mathbb{X}_i}$ that passed the previous screening stage $S_{i-1}$, where $\mathbb{X}_i$ is given by:
\begin{equation}
    \mathbb{X}_{i} = \left\{x \mid  x \in \mathbb{X}_{i-1} \ \mbox{and} \ f_{i-1}(x) \ge \lambda_{i-1}  \right\}.
\end{equation}
By definition, we have $\mathbb{X}_{i} \triangleq \mathbb{X}$, which contains the entire set of molecules to be screened.
At stage ${S_i}$, every molecule $x \in \mathbb{X}_{i}$ whose property score ${y_i = f_i\left(x\right)}$ is below the threshold ${\lambda_i}$ is discarded such that only the remaining molecules $x \in \mathbb{X}_{i+1}$ that meet or exceed this threshold are passed on to the next stage ${S_{i+1}}$.
We assume that all molecules in $\mathbb{X}_i$ at each stage ${S_i}$ are batch-processed to select the set of molecules $\mathbb{X}_{i+1}$ that will be passed to the subsequent stage $S_{i+1}$, as it is often done in practice~\cite{wilmer2012large,gupta2020computational,kim2021universal}.

Although every stage $S_i$ in the screening pipeline performs a down-selection of the molecules by assessing their molecular property based on the surrogate model $f_i(x)$ and comparing it against the threshold $\lambda_i$, we assume only the threshold values $\lambda_1, \cdots, \lambda_{i-1}$ of the first $N-1$ stages will need to be determined while the threshold ${\lambda_N}$ for the last screening stage ${S_N}$ is predetermined. This reflects how such screening pipelines are utilized in real-world scenarios. For example, in the IMPECCABLE pipeline~\cite{saadi2020impeccable}, as well as in many other computational drug discovery pipelines, potentially effective lead compounds that pass the earlier stages based on efficient but less accurate models will be assessed using computationally expensive yet highly accurate molecular dynamics (MD) simulations to evaluate the binding affinity against the target. Only the molecules whose binding affinity estimated by the MD simulations exceeds a reasonably high threshold set by domain experts may be further assessed experimentally, in order not to unnecessarily waste the available resources, considering that such experimental validation is typically costly, time-consuming, and labor-intensive. Similarly, in a materials screening pipeline, the last screening stage may involve expensive calculations based on density-functional theory (DFT), a quantum mechanical modeling scheme that is widely used for predicting material properties~\cite{liu2015high,kim2016first,kim2016thermodynamic,liu2017self,park2017systematic,kang2018density,sood2018electrochemical,sood2018electrochemical2,zhu2018boron}.

Based on this setting, our primary objective is to predict the optimal screening policy ${\boldsymbol{\psi}^\ast  = \left[\lambda_1^\ast, \lambda_2^\ast, \dots, \lambda_{N-1}^\ast \right]}$ that leads to the optimal operation of the HTVS pipeline. We consider two different scenarios. In the first scenario, we assume that the total computational budget for screening the candidate molecules is fixed, where the design goal would then be to identify the optimal screening policy that maximizes the screening throughput, namely, the percentage (or number) of potential molecules that meet or exceed the qualification in the last stage $S_N$ (\textit{i.e.}, $f_N(x)\ge \lambda_N$). In the second scenario, we consider the case when the computational budget is not fixed and where the goal is to design the optimal policy that simultaneously maximizes the throughput while minimizing the overall computational cost.

Figure~\ref{fig1} (right) shows a flowchart summarizing the proposed approach for identifying the optimal screening policy ${\boldsymbol{\psi}^\ast  = \left[\lambda_1^\ast, \lambda_2^\ast, \dots, \lambda_{N-1}^\ast \right]}$ for the optimal operation of a given HTVS pipeline under the two screening scenarios described above. First, we estimate the joint distribution ${p \left(y_1, y_2, \dots, y_N \right)}$ of the predictive scores from the $N$ stages based on the available training data. In case the probability density function (PDF) ${p \left(y_1, y_2, \dots, y_N \right)}$ is known \textit{a priori}, this PDF estimation step will not be required.
Given ${p \left(y_1, y_2, \dots, y_N \right)}$, we can predict the optimal screening policy ${ \boldsymbol{\psi}^\ast  = \left[\lambda_1^\ast, \lambda_2^\ast, \dots, \lambda_{N-1}^\ast \right]}$ that leads to the optimal operational performance of the HTVS pipeline. Specifically, in case the total computational budget ${C}$ is fixed, we find the optimal policy ${\boldsymbol{\psi}^\ast  = \left[\lambda_1^\ast, \lambda_2^\ast, \dots, \lambda_{N-1}^\ast \right]}$ that maximizes the screening throughput of the pipeline--\textit{i.e.}, the proportion of molecules that pass the last (and the most stringent/accurate) screening stage that meet the condition ${f_N \left(x \right) \ge \lambda_N}$--under the budget constraint ${C}$. The formal definition of the given optimization problem is shown in Eq.~\eqref{eq:optimization}. Otherwise, we predict optimal screening policy ${ \boldsymbol{\psi}^\ast  = \left[\lambda_1^\ast, \lambda_2^\ast, \dots, \lambda_{N-1}^\ast \right]}$ that jointly optimizes the throughput and computational resource based on a weighted objective function that balances the throughput and the computational cost.
The formal definition of this joint optimization problem can be found in Eq.~\eqref{eq:optimization:joint}.
In this case, the balancing weight ${\alpha}$ can be used to trade throughput for computational efficiency, or vice versa. We note that the training dataset is only used for estimating the PDF ${p \left(y_1, y_2, \dots, y_N \right)}$ and not (directly) for finding the optimal screening policy. In fact, the optimal policy ${\boldsymbol{\psi}^\ast  = \left[\lambda_1^\ast, \lambda_2^\ast, \dots, \lambda_{N-1}^\ast \right]}$ is determined by a function of up to three parameters: the joint score distribution ${p \left(y_1, y_2, \dots, y_N \right)}$, ${\lvert \mathbb{X} \rvert}$ (the number of potential molecules to be screened), and the total computational budget ${C}$ (in the first screening scenario, where the computational budget is assumed to be limited).

\subsection*{Learning the output relations across multiple stages in a high-throughput virtual screening pipeline}
\label{sec:learningthejointscoredistribution}
As shown in Fig.~\ref{fig1}, the proposed optimization framework that identifies the optimal screening policy ${ \boldsymbol{\psi}^\ast}$ takes a two-phase approach. In the first phase, we estimate the joint score distribution ${p\left(y_1, y_2, \dots, y_N \right)}$. Based on the estimated score distribution, we find the optimal screening policy ${ \boldsymbol{\psi}^\ast}$ that maximizes the screening performance. To ensure good screening performance, accurate estimation of the joint score distribution ${p\left(y_1, y_2, \dots, y_N \right)}$ is crucial. In this study, we performed a spectral estimation under the assumption that the joint score distribution follows a multivariate Gaussian mixture model and estimated the parameters via the expectation-maximization (EM) scheme~\cite{dempster1977maximum}.


\subsubsection*{Finding the optimal screening policy under computational budget constraint}
\label{firstFramework}
A formal objective is to compute screening threshold value ${\lambda_i}$ at stage ${S_i}$ for ${i = 1, 2, \dots, N-1}$, such that the total number of the detected potential candidates in the set ${\mathbb{Y}}$, where the candidates meet the target criteria based on the score in the last stage ${S_N}$ and the pre-specified screening threshold ${\lambda_N}$, is maximized under a given computational budget ${C}$. The relationship among the predictive scores from all stages ${S_1, S_2, \dots, S_N}$ is captured by their joint score distribution, which is denoted as ${p\left(y_1, y_2, \dots, y_N \right)}$. Based on this joint score distribution, we define the following reward function ${r\left(\boldsymbol{\lambda}\right)}$ according to policy ${\boldsymbol{\lambda} = \left[\lambda_1, \lambda_2, \dots, \lambda_{N} \right]}$ of the stages ${S_i}$, ${i = 1, 2, \dots, N}$, as follows:
\begin{equation}
    r\left(\boldsymbol{\lambda}\right) = \idotsint\limits_{[\lambda_{N}, \lambda_{N-1}, \dots, \lambda_1]}^{\boldsymbol{\infty}} p\left(y_1,y_2,\dots,y_{N}\right) \,dy_1 dy_2 \cdots dy_{N}.
\end{equation}
%
We can find the optimal screening policy ${ \boldsymbol{\psi}^\ast  = \left[\lambda_1^\ast, \lambda_2^\ast, \dots, \lambda_{N-1}^\ast \right]}$ to be applied to the first $N-1$ stages (${S_i}$, ${i = 1, 2, \dots, N-1}$) that maximizes the reward ${\lvert \mathbb{Y} \rvert}$ by solving the constrained optimization problem shown below:
\begin{align}
\boldsymbol{\psi}^\ast = \argmax_{\boldsymbol{\psi} \in \mathbb{R}^{N-1}}  & \quad r\left(\left[\boldsymbol{\psi}, \lambda_N \right]\right) \label{eq:optimization}\\
\textrm{s.t.} &\quad \sum_{i=1}^{N} c_i \lvert \mathbb{X}_i\rvert\  \le C, \nonumber
\end{align}
where ${ \lvert \mathbb{X}_i \rvert }$ is the number of molecules that passed the previous stages from ${S_1}$ to ${S_{i-1}}$. Formally, $\lvert \mathbb{X}_i \rvert$ is defined as:
\begin{equation}
\lvert \mathbb{X}_i \rvert =  \lvert \mathbb{X} \rvert\idotsint\limits_{[\lambda_{i-1}, \lambda_{i-2},\dots, \lambda_1]}^{\boldsymbol{\infty}} p_{1:i-1}\left(y_1,y_2,\dots,y_{i-1}\right) \,dy_1 dy_2 \cdots dy_{i-1},
\end{equation}
where $p_{1:i-1}$ denotes the marginal score distribution for $y_1, \cdots, y_{i-1}$, which can be obtained by marginalizing $p(\cdot)$ over ${y_i}$ to ${y_N}$.


\subsubsection*{Joint optimization of the screening policy for throughput and computational efficiency}
\label{secondFramework}

In many real-world screening problems, including drug or material screening, the total computational budget for screening may not be fixed, and one may want to jointly optimize for both screening throughput as well as computational efficiency of screening. In such scenarios, we can formulate a joint optimization problem to find the best screening policy that strikes the optimal balance between throughput and efficiency:
\begin{equation}
\boldsymbol{\psi}^\ast = \argmin_{\boldsymbol{\psi} \in \mathbb{R}^{N-1}} \alpha \bar{r}\left(\left[\boldsymbol{\psi},  \lambda_N\right]\right) + \left(1-\alpha\right) \bar{h}\left(\left[ \boldsymbol{\psi},\lambda_N\right]\right).  \label{eq:optimization:joint}
\end{equation}
The weight parameter ${\alpha \in \left[0, 1\right]}$ determines the relative importance between the relative reward function ${\bar{r}\left(\left[\boldsymbol{\psi}, \lambda_N\right]\right)}$ and the normalized total cost function ${\bar{h}\left(\left[ \boldsymbol{\psi},\lambda_N\right]\right)}$ defined as follows:
\begin{align}
\bar{r} \left(\left[\boldsymbol{\psi}, \lambda_N\right]\right) &= \frac{r\left(\left[-\boldsymbol{\infty}, \lambda_N \right]\right) - r\left(\left[\boldsymbol{\psi}, \lambda_N \right]\right)}{r\left(\left[-\boldsymbol{\infty}, \lambda_N \right]\right)}\\
&= \frac{\int_{\lambda_N}^{\infty} p_N\left(y_N \right)\,dy_{N} - r\left(\left[\boldsymbol{\psi}, \lambda_N \right]\right)}{\int_{\lambda_N}^{\infty} p_N\left(y_N \right)\,dy_{N}},
\end{align}
\begin{equation}
\bar{h}\left(\left[\boldsymbol{\psi}, \lambda_N\right]\right) = \frac{1}{N \lvert \mathbb{X} \rvert \max_{i}c_i} \sum_{i=1}^{N} c_i\lvert \mathbb{X}_i \rvert.
\end{equation}
Note that $p_N$ is the marginal score distribution for $y_N$, which is obtained by marginalizing $p(\cdot)$ over ${y_1}$ to ${y_{N-1}}$.

\section{Results}
\label{sec:results}
We validate the proposed optimization framework based on both synthetic and real data. First, we evaluate the performance of our optimization framework based on a four-stage HTVS pipeline, where the joint probability distribution of the predictive scores is assumed to be known. Next, we construct an HTVS pipeline for lncRNAs by interconnecting existing lncRNA prediction algorithms with different prediction accuracy and computational complexity. In this example, the joint distribution of the predictive scores from the different algorithms at different stages is learned from training data, based on which the proposed HTVS optimal framework is used to identify the optimal screening policy.

\subsection*{Comprehensive performance analysis of the HTVS pipeline optimization framework}

For comprehensive performance analysis of the proposed HTVS pipeline optimization framework, we consider a synthetic HTVS pipeline with $N=4$ stages, where the joint PDF of the predictive scores from all stages is assumed to be known. We vary the correlation levels between the scores from neighboring stages to investigate the overall impact on the performance of the optimized HTVS pipeline. 

Specifically, we assume that the computational cost for screening a single molecule is ${1}$ at stage $S_1$, ${10}$ at $S_2$, ${100}$ at $S_3$, and ${1,000}$ at $S_N$. As the per-molecule screening cost is fairly different across stages, the given setting for the synthetic HTVS pipeline allows us to clearly see the impact and significance of optimal decision-making on the overall throughput and accuracy of the screening pipeline.

Here we consider the case when we have complete knowledge of the joint score distribution ${p \left(y_1, y_2, y_3, y_4 \right)}$. The score distribution is assumed to be a multivariate uni-modal Gaussian distribution ${\mathcal{G}\left(\mathbf{0}, \mathbf{\Sigma}\left(\rho\right)\right)}$, where the covariance matrix ${\mathbf{\Sigma}\left(\rho\right)}$ is a Toeplitz matrix defined as follows:
\begin{equation}
    \mathbf{\Sigma}\left(\rho\right) = 
    \begin{bmatrix}
        1       & \rho     & \rho-0.1 & \rho-0.2 \\
        \rho       & 1     & \rho     & \rho-0.1 \\
        \rho-0.1   & \rho     & 1     & \rho     \\
        \rho-0.2   & \rho-0.1 & \rho     & 1
    \end{bmatrix},
\end{equation}
where $\rho$ is the correlation between neighboring stages $S_{i}$ and $S_{i+1}$ for $i=1,2,3$. We assumed that the score correlation is lower between stages that are further apart, which is typically the case in real screening pipelines that consist of multi-fidelity models.

The primary objective of the HTVS pipeline is to maximize the number of potential candidates that satisfy the final screening criterion (\textit{i.e.}, ${f_4\left(x\right) \ge \lambda_4}$) based on the highest fidelity model at stage $S_4$ while minimizing the total computational cost induced by the entire screening pipeline. The total number of all candidate molecules in the initial set $\mathbb{X}$ is assumed to be ${10^5}$. We assume that we are given ${\lambda_4 = 3.0902}$ as prior information set by a domain expert, which results in ${100}$ molecules (among ${10^5}$ in $\mathbb{X}$) that satisfy the final screening criterion. We validate the proposed HTVS optimization framework for two cases: first, for $\rho=0.8$, where the neighboring stages yield scores that are highly correlated, and next, for $\rho=0.5$ where the correlation is relatively low. Performance analysis results based on various other covariance matrices can be found in the supplementary materials. 

\subsubsection*{Performance of the optimized HTVS pipeline under computational budget constraint}

\begin{figure}[h!]
\includegraphics[width=\textwidth]{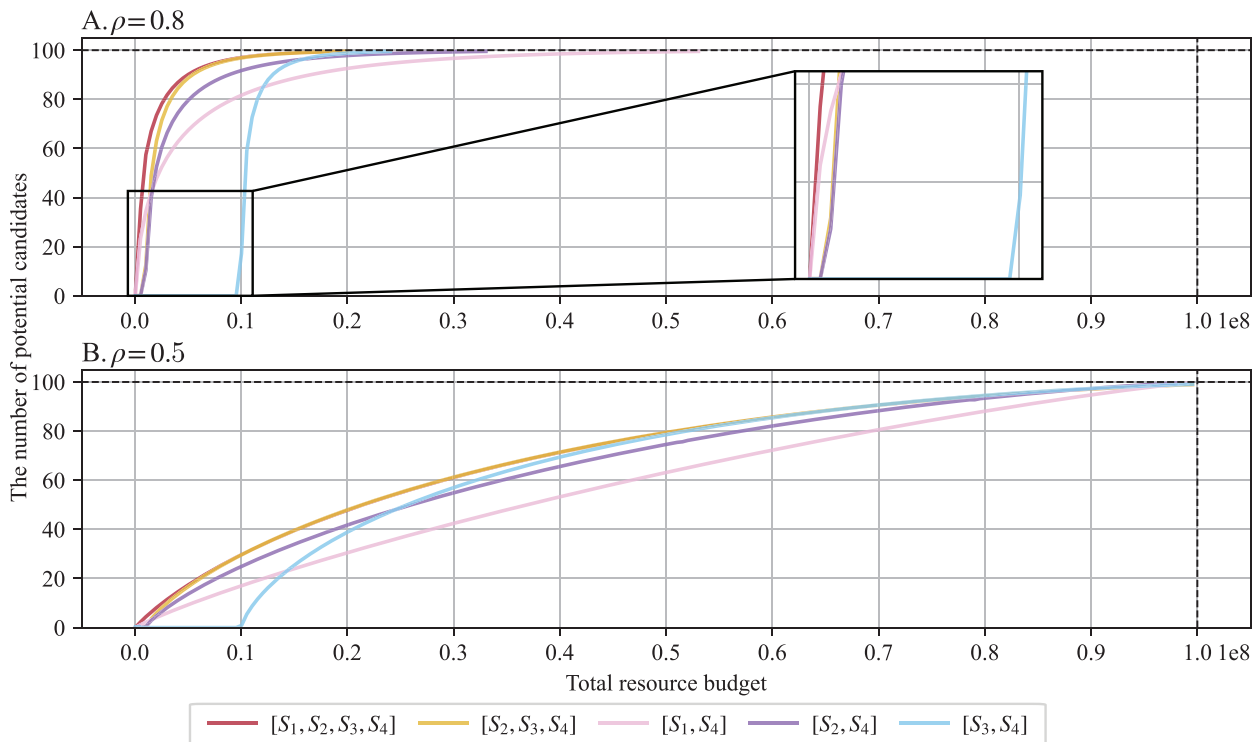}
\caption{Performance assessment of the optimized HTVS pipelines. The number of candidate molecules that meet the desired screening criterion is shown as a function of the available computational budget. Results are shown for the case when the stages are highly correlated (A, ${\rho = 0.8}$) as well as when they have relatively low correlation (B. ${\rho = 0.5}$). The performance of the best-performing $4$-stage pipeline and the best-performing $3$-stage pipeline is shown. For comparison, we also show the performance of all $2$-stage pipelines. Note that only the best-performing configurations are shown for ${N \ge 3}$.}
\label{fig2}
\end{figure}

Figure~\ref{fig2} shows the performance evaluation results for different HTVS pipeline structures optimized via the proposed framework under a fixed computational resource budget. The total number of the desirable candidates detected by the pipeline is shown as a function of the available computational budget for two cases: (A) HTVS pipeline that consists of highly-correlated stages (\textit{i.e.}, ${\rho = 0.8}$) and (B) HTVS pipeline comprised of stages with lower correlation (\textit{i.e.}, ${\rho = 0.5}$). The black horizontal and vertical dashed lines depict the total number of true candidates that meet the screening criterion (${100}$ in this simulation) and the total computational budget required when screening all molecules in $\mathbb{X}$ only based on the last stage ${S_4}$ (\textit{i.e.}, the highest-fidelity and most computationally expensive model), respectively. Figure~\ref{fig2} shows the performance of the best-performing $N=4$ stage pipeline and that of the best-performing $N=3$ pipeline. Additionally, the performance of all $N=2$ stage pipelines is shown for comparison.

First, as shown in Fig.~\ref{fig2}.A, the performance curves of the pipelines consisting of only two stages (shown in pink, purple, and light blue lines) demonstrate how each of the lower-fidelity stages $S_1$--$S_3$ improves the screening performance when combined with the highest-fidelity stage $S_4$ and performance-optimized by our proposed framework.
As shown in Fig.~\ref{fig2}.A, the correlation between the lower-fidelity/lower-complexity stage ${S_i}$, ${i = 1, 2, \dots, N- 1}$, at the beginning of the HTVS pipeline and highest-fidelity/highest-complexity stage ${S_N}$ at the end of the pipeline has a significant impact on the slope of the performance curve. For example, in the two-stage pipeline ${\left[S_3, S_4 \right]}$, where the two stages are highly correlated to each other, we can observe the steepest performance improvement as the available computational budget increases.
On the other hand, for the two-stage pipeline ${\left[S_1, S_4 \right]}$ which consists of less correlated stages, the performance improvement is more moderate in comparison as the available computational budget increases. 
Note that the minimum required computational budget to screen all candidates is larger for the pipeline ${\left[S_3, S_4 \right]}$ compared to that for ${\left[S_1, S_4 \right]}$, which is due to the assumption that all candidates are batch-processed at each stage. For example, with the minimum budget needed by pipeline ${\left[S_3, S_4 \right]}$ to screen all candidates, the other pipelines ${\left[S_1, S_4 \right]}$ and ${\left[S_2, S_4 \right]}$ are capable of completing the screening and detecting more than ${80\%}$ of the desirable candidates. Nevertheless, the detection performance improves with the increasing computational budget for all two-stage pipelines.

It is important to note that we can in fact simultaneously attain the advantage of using a lower-complexity stage (\textit{e.g.}, ${\left[S_1, S_4 \right]}$) that allows a ``quick-start'' with a small budget as well as the merit of using a higher-complexity stage (\textit{e.g.}, ${\left[S_3, S_4 \right]}$) for rapid performance improvement with the budget increase by constructing a multi-stage HTVS pipeline and optimally allocating the computational resources according to our proposed optimization framework.
This can be clearly seen in the performance curve for the four-stage pipeline ${\left[S_1, S_2, S_3, S_4\right]}$ (shown in the red solid line). The optimized four-stage pipeline consistently outperforms all other pipelines across all budget levels. 
Specifically, the optimized pipeline ${\left[S_1, S_2, S_3, S_4\right]}$ quickly evaluated all the molecular candidates in $\mathbb{X}$ through the most efficient stage ${S_1}$ and sharply improved the screening performance through the utilization of more complex yet also more accurate subsequent stages in the HTVS pipeline in a resource-optimized manner. For example, the optimized four-stage pipeline detected ${97\%}$ of the desirable candidates that meet the target criterion at only ${10\%}$ of the total computational cost that would be required if one used only the last stage (which we refer to as the ``original cost''). To detect ${99\%}$ of the desired candidates, the optimized four-stage pipeline ${\left[S_1, S_2, S_3, S_4\right]}$ would need only about ${14\%}$ of the original cost. 

Among all three-stage pipelines (\textit{i.e.}, ${N = 3}$), pipeline ${\left[S_2, S_3, S_4\right]}$ yielded the best performance when performance-optimized using our proposed optimization framework (orange solid line in Fig.~\ref{fig2}.A). As we can see in Fig.~\ref{fig2}.A, the screening performance sharply increases as the available computational budget increases, thanks to the high correlation between $S_4$ and the prior stages ${S_2}$ and ${S_3}$. However, due to the higher computational complexity of ${S_2}$ compared to that of ${S_1}$, the optimized pipeline ${\left[S_2, S_3, S_4\right]}$ required a higher minimum computational budget for screening all candidate molecules compared to the minimum budget needed by a pipeline that begins with ${S_1}$. Despite this fact, when the first stage ${S_2}$ in this three-stage HTVS pipeline is replaced by the more efficient ${S_1}$, our simulation results (see Fig.~S41 in the supplementary materials) show that the screening performance improves relatively moderately as the budget increases. 
Empirically, when all stages are relatively highly correlated to each other, the best strategy for constructing the HTVS pipeline appears to place the stages in increasing order of complexity and optimally allocate the computational resources to maximize the return-on-computational-investment (ROCI). In fact, this observation is fairly intuitive and also in agreement with how screening pipelines are typically constructed in real-world applications. 

Figure~\ref{fig2}.B shows the performance evaluation results of the HTVS pipelines, where the screening stages are moderately correlated to each other (${\rho = 0.5}$). Results are shown for different pipeline configurations, where the screening policy is optimized using the proposed framework to maximize the ROCI. Overall, the performance trends were nearly identical to those shown in Fig.~\ref{fig2}.A, although the overall performance is lower compared to the high correlation scenario (${\rho = 0.8}$) as expected. 
While the screening performance of the optimized HTVS pipeline is not as good as the high-correlation scenario, the multi-stage HTVS pipeline with the optimized screening policy still provides a much better trade-off between the computational cost for screening and the detection performance. For example, if we were to use only the highest-fidelity model in $S_4$ for screening, the only way to trade accuracy for reduced resource requirements would be to randomly sample the candidate molecules from $\mathbb{X}$ and screen the selected candidates. The performance curve in this case would be a straight line connecting $\left(0,0\right)$ and $\left(10^8,100\right)$, below most of the performance curves for the optimized pipeline approach shown in Fig.~\ref{fig2}.B. As in the previous case ($\rho=0.8$), the best pipeline configuration was to interconnect all four stages, where the stages are connected to each other in increasing order of complexity.

\subsubsection*{Performance of the HTVS pipeline jointly optimized for throughput and computational efficiency}

Table~\ref{tab1} shows the performance of the various HTVS pipeline configurations, where the screening policy was jointly optimized for both throughput and computational efficiency. The joint optimization problem is formally defined in Eq.~\eqref{eq:optimization:joint}, and $\alpha$ was set to $0.5$ in these simulations.
As a reference, the first row (configuration ${\left[S_4\right]}$) shows the performance of solely relying on the last stage ${S_4}$ for screening the molecules without utilizing a multi-stage pipeline. The effective cost is defined as the total computational cost divided by the total number of molecules detected by the screening pipeline that satisfy the target criterion (\textit{i.e.}, average computational cost per detected candidate molecule).
The computational savings of a given pipeline configuration is calculated by comparing its effective cost to that of the reference configuration (\textit{i.e.}, ${\left[S_4\right]}$).
As we can see in Table~\ref{tab1}, our proposed HTVS pipeline optimization framework was able to significantly improve the overall screening performance across all pipeline configurations in a highly robust manner.
For example, for $\rho=0.8$, the optimized pipelines consistently led to computational savings ranging from ${76.20\%}$ to ${86.64\%}$ compared to the reference, while detecting ${94\sim99 \%}$ of the desired candidates that meet the target criterion.
Although the overall efficiency of the HTVS pipelines slightly decreases when the neighboring stages are less correlated (${\rho=0.5}$), the pipelines were nevertheless effective in saving computational resources. As shown in Table~\ref{tab1}, the optimized HTVS pipelines detected ${89\% \sim 96\%}$ of all desired candidate molecules with computational savings ranging between ${36.46\%}$ and ${54.87\%}$.

\begin{table*}
\scriptsize
\centering
\begin{tabular}{ccccccccccc}
\hline
\multirow{2}{*}{Configuration}  & & \multicolumn{4}{c}{High correlation ${\left(\rho=0.8\right)}$}&&\multicolumn{4}{c}{Low correlation ${\left(\rho=0.5\right)}$}\\ 
\cline{3-6} \cline{8-11}& & \begin{tabular}[c]{@{}c@{}}Potential\\candidates\end{tabular} & Total cost & \begin{tabular}[c]{@{}c@{}}Effective\\cost\end{tabular} & \begin{tabular}[c]{@{}c@{}}Comp.\\savings\end{tabular} &  & \begin{tabular}[c]{@{}c@{}}Potential\\candidates\end{tabular} & Total cost & \begin{tabular}[c]{@{}c@{}}Effective\\cost\end{tabular} & \begin{tabular}[c]{@{}c@{}}Comp.\\savings\end{tabular} \\ \hline
{${\left[S_4\right]}$}                  & & $100$       & ${100,000,000}$   & ${1,000,000}$ & ${0\%}$       & & ${100}$     & ${100,000,000}$   & ${1,000,000}$ & ${0\%}$ \\ \hline
{${\left[S_1, S_4\right]}$}             & & {${94}$}    & ${22,372,654}$    & ${238,007}$   & ${76.20\%}$   & & {${89}$}    & {${56,551,129}$}  & ${635,406}$   & ${36.46\%}$ \\
{${\left[S_2, S_4\right]}$}             & & {${96}$}    & ${15,511,702}$    & ${161,580}$   & ${83.84\%}$   & & {${90}$}    & {${43,620,751}$}  & ${484,675}$   & ${51.53\%}$ \\
{${\left[S_3, S_4\right]}$}             & & {${98}$}    & ${18,152,330}$    & ${185,228}$   & ${81.48\%}$   & & {${92}$}    & {${41,522,035}$}  & ${451,326}$   & ${54.87\%}$ \\
{${\left[S_1, S_2, S_4\right]}$}        & & {${97}$}    & ${17,890,176}$    & ${184,435}$   & ${81.56\%}$   & & {${94}$}    & ${53,340,817}$    & ${567,456}$   & ${43.25\%}$ \\
{${\left[S_1, S_3, S_4\right]}$}        & & {${98}$}    & ${14,451,644}$    & ${147,466}$   & ${85.25\%}$   & & {${94}$}    & ${47,550,232}$    & ${505,854}$   & ${49.41\%}$ \\
{${\left[S_2, S_1, S_4\right]}$}        & & {${97}$}    & ${18,291,054}$    & ${188,568}$   & ${81.14\%}$   & & {${94}$}    & ${53,513,582}$    & ${569,293}$   & ${43.07\%}$ \\
{${\left[S_2, S_3, S_4\right]}$}        & & {${98}$}    & ${13,089,779}$    & ${133,569}$   & ${86.64\%}$   & & {${94}$}    & ${44,534,328}$    & ${473,769}$   & ${52.62\%}$ \\
{${\left[S_3, S_1, S_4\right]}$}        & & {${99}$}    & ${19,505,326}$    & ${197,023}$   & ${80.30\%}$   & & {${94}$}    & ${48,708,112}$    & ${518,171}$   & ${48.18\%}$ \\
{${\left[S_3, S_2, S_4\right]}$}        & & {${99}$}    & ${19,522,312}$    & ${197,195}$   & ${80.28\%}$   & & {${94}$}    & ${47,966,605}$    & ${510,283}$   & ${48.97\%}$ \\
{${\left[S_1, S_2, S_3, S_4\right]}$}   & & {${99}$}    & ${14,147,264}$    & ${142,902}$   & ${85.71\%}$   & & {${96}$}    & ${50,336,621}$    & ${524,340}$   & ${47.57\%}$ \\
{${\left[S_1, S_3, S_2, S_4\right]}$}   & & {${99}$}    & ${15,939,108}$    & ${161,001}$   & ${83.90\%}$   & & {${96}$}    & ${52,704,450}$    & ${549,005}$   & ${45.10\%}$ \\
{${\left[S_2, S_1, S_3, S_4\right]}$}   & & {${99}$}    & ${14,348,794}$    & ${144,937}$   & ${85.51\%}$   & & {${96}$}    & ${50,366,503}$    & ${524,651}$   & ${47.53\%}$ \\
{${\left[S_2, S_3, S_1, S_4\right]}$}   & & {${99}$}    & ${14,335,230}$    & ${144,800}$   & ${85.52\%}$   & & {${96}$}    & ${50,411,458}$    & ${525,119}$   & ${47.49\%}$ \\
{${\left[S_3, S_1, S_2, S_4\right]}$}   & & {${99}$}    & ${20,560,571}$    & ${207,682}$   & ${79.23\%}$   & & {${96}$}    & ${53,249,970}$    & ${554,687}$   & ${44.53\%}$ \\
{${\left[S_3, S_2, S_1, S_4\right]}$}   & & {${99}$}    & ${20,560,299}$    & ${207,680}$   & ${79.23\%}$   & & {${96}$}    & ${53,215,674}$    & ${554,330}$   & ${44.57\%}$ \\ \hline 
\end{tabular}
\caption{Performance comparison of various high-throughput virtual screening (HTVS) pipeline structures jointly optimized via the proposed framework (${\alpha = 0.5}$).}
\label{tab1}
\end{table*}

For further evaluation of the proposed framework, we performed additional experiments based on the four-stage pipeline ${\left[S_1, S_2, S_3, S_4\right]}$. In this experiment, we first investigated the impact of $\alpha$ on the screening performance. Next, we compared the performance of the optimal screening policy with the performance of a baseline policy that mimics a typical screening scenario in real-world applications (\textit{e.g.}, see~\cite{saadi2020impeccable}). The baseline policy selects the top ${R_s}\%$ candidate molecules at each stage and passes them to the next stage while discarding the rest. This baseline screening policy is agnostic of the joint score distribution of the multiple stages in the HTVS and aims to reduce the overall computational cost by passing only the top candidates to subsequent stages that are more costly. Similar strategies are in fact often adopted in practice due to their simplicity. In our simulations, we assumed the proportion ${R_s}$ is uniform across the screening stages.
The performance evaluation results are summarized in Table~\ref{tab2}. When the neighboring stages were highly correlated (${\rho = 0.8}$), the optimized pipelines detected ${100}$, ${99}$, and ${96}$ candidate molecules at a total cost of ${19,727,704}$; ${14,147,264}$; and ${10,926,901}$, respectively. Interestingly, when ${\alpha}$ was reduced from ${0.75}$ to ${0.25}$ (\textit{i.e.}, trading accuracy for higher efficiency), the number of detected candidate molecules decreased only by $4$ (\textit{i.e.}, from $100$ to $96$), while leading to an additional computational savings of $8$ percentage points (\textit{i.e.}, from $80.27\%$ to $88.62\%$). 
On the other hand, the performance of the baseline screening policy was highly unpredictable and very sensitive to the choice of ${R_s}$. 
For example, although the baseline with ${R_s=0.75}$ found all the potential candidates, the effective cost of the baseline was significantly higher than that of the proposed optimized pipeline with $[\alpha =0.75]$. For ${R_s = 0.5}$, the baseline detected ${98}$ potential candidates (out of $100$) with a total cost of ${15,599,934}$, which was higher than the total cost of the optimized pipeline that detected ${99}$ potential candidates. The baseline pipelines with ${R_s =0.1}$ and ${0.25}$ selected ${26\%}$ and ${78\%}$ of the potential candidates, respectively. Considering that the primary goal of such a pipeline is to detect the largest number of potential candidates in a computationally efficient manner, these results clearly show that this baseline screening scheme that mimics conventional screening pipelines resulted in unreliable and suboptimal performance even when the neighboring stages were highly correlated to each other.
While the baseline may lead to reasonably good performance for certain ${R_s}$, it is important to note that we cannot determine the optimal ${R_s}$ in advance as the approach is agnostic to the relationships between different stages. As a result, the application of this baseline screening pipeline may significantly degrade the screening performance in practice. When the correlation between the neighboring stages was relatively low (${\rho=0.5}$), the overall performance of the proposed pipeline degraded as expected. In this case, the pipeline jointly optimized for screening accuracy as well as efficiency with $\alpha$ set to ${0.75}$, ${0.5}$, and ${0.25}$ detected ${99}$, ${96}$, and ${86}$ potential candidates with the computational cost of ${71,836,915}$, ${50,336,621}$, and ${28,563,886}$, respectively. As in the high correlation case, the performance of the baseline scheme significantly varied and was sensitive to the choice of ${R_s}$.

\begin{table*}
\scriptsize
\centering
\begin{tabular}{ccccccccccc}
\hline
\multirow{2}{*}{approach}  & & \multicolumn{4}{c}{High correlation ${\left(\rho=0.8\right)}$}&&\multicolumn{4}{c}{Low correlation ${\left(\rho=0.5\right)}$}\\ 
\cline{3-6} \cline{8-11}& & \begin{tabular}[c]{@{}c@{}}Potential\\candidates\end{tabular} & Total cost & \begin{tabular}[c]{@{}c@{}}Effective\\cost\end{tabular} & \begin{tabular}[c]{@{}c@{}}Comp.\\savings\end{tabular} &  & \begin{tabular}[c]{@{}c@{}}Potential\\candidates\end{tabular} & Total cost & \begin{tabular}[c]{@{}c@{}}Effective\\cost\end{tabular} & \begin{tabular}[c]{@{}c@{}}Comp.\\savings\end{tabular} \\ \hline
Proposed (${\alpha = 0.75}$)            & & {${100}$}    & ${19,727,704}$    & ${197,277}$   & ${80.27\%}$   & & {${99}$}    & ${71,836,915}$    & ${725,625}$   & ${27.44\%}$\\
Proposed (${\alpha = 0.5}$)             & & {${99}$}    & ${14,147,264}$    & ${142,902}$   & ${85.71\%}$   & & {${96}$}    & ${50,336,621}$    & ${524,340}$   & ${47.57\%}$\\
Proposed (${\alpha = 0.25}$)            & & {${96}$}    & ${10,926,901}$    & ${113,822}$   & ${88.62\%}$   & & {${86}$}    & ${28,563,886}$    & ${332,138}$   & ${66.79\%}$\\\hline
Baseline ($R_s = 75\%$)                 & & {${100}$}   & ${48,966,384}$    & ${489,664}$   & ${51.03\%}$   & & {${93}$}    & {${48,662,387}$}  & ${523,251}$   & ${47.67\%}$ \\
Baseline ($R_s = 50\%$)                 & & {${98}$}    & ${15,599,934}$    & ${159,183}$   & ${84.08\%}$   & & {${69}$}    & {${15,600,165}$}  & ${226,089}$   & ${77.39\%}$ \\
Baseline ($R_s = 25\%$)                 & & {${78}$}    & ${2,537,498}$     & ${32,532}$    & ${96.75\%}$   & & {${28}$}    & {${2,537,516}$}   & ${90,626}$    & ${90.94\%}$ \\
Baseline ($R_s = 10\%$)                 & & {${26}$}    & ${400,000}$       & ${15,385}$    & ${98.46\%}$   & & {${6}$}     & {${400,008}$}     & ${66,668}$    & ${93.33\%}$ \\
\hline 
\end{tabular}
\caption{Performance comparison between the proposed pipeline {${\left[S_1, S_2, S_3, S_4\right]}$} jointly optimized for throughput and computational efficiency (with various ${\alpha}$) and the baseline pipeline (with different screening ratio ${R_s}$) in terms of the total number of detected potential candidates after screening and the computational cost induced.}
\label{tab2}
\end{table*}


\subsection*{Performance evaluation of the optimized HTVS pipeline for screening long non-coding RNAs}

To demonstrate the efficacy of the proposed optimization framework in a real-world application, we considered an optimal computational screening campaign for the identification of long non-coding RNAs (lncRNAs). Given a large number of RNA transcripts, the goal is to efficiently and accurately detect lncRNA transcripts through an HTVS pipeline. In recent years, interests in lncRNAs have been constantly increasing in relevant research communities, as there is growing evidence that lncRNAs and their roles in various biological processes are closely associated with the development of complex and often hard-to-treat diseases including Alzheimer's diseases~\cite{ng2013long,tan2013non,luo2016long}, cardiovascular diseases~\cite{congrains2012genetic,xue2016g}, as well as several types of cancer~\cite{yang2014lncrna,shi2015critical,peng2017lncrna,carlevaro2020cancer}.
RNA sequencing techniques are nowadays routinely used to investigate the main functional molecules and their molecular interactions responsible for the initiation, progression, and manifestation of such complex diseases. Consequently, the accurate detection of lncRNA transcripts from a potentially huge number of sequenced RNA transcripts is a fundamental step in studying lncRNA-disease association.
While several lncRNA prediction algorithms have been developed so far~\cite{wang2013cpat,li2014plek,kang2017CPC2,han2019lncfinder}, each of which with its own pros and cons, no HTVS pipeline has been proposed to date for fast and reliable screening of lncRNAs.
In this study, we first constructed a lncRNA HTVS pipeline by interconnecting four existing prediction algorithms--CPC$2$ (Coding Potential Calculator $2$)~\cite{kang2017CPC2}, CPAT (Coding Potential Assessment Tool)~\cite{wang2013cpat}, PLEK (Predictor of LncRNAs and mEssenger RNAs based on an improved $k$-mer scheme)~\cite{li2014plek}, and LncFinder~\cite{han2019lncfinder}. Next, we estimated the joint probability distribution of the predictive scores from the given algorithms. The estimated score distribution was then \tcr{used} to derive the optimal screening policy according to our proposed HTVS optimization framework. In we describe each of these steps in further detail and present the performance evaluation results based on real RNA transcripts in the GENCODE database~\cite{frankish2021gencode}.

\subsubsection*{Dataset and preprocessing}

We collected the nucleotide sequences of \textit{Homo sapiens} RNA transcripts from GENCODE v$38$ (May $5$, $2021$)~\cite{frankish2021gencode}, which consists of $48,752$ lncRNA sequences and $106,143$ protein-coding sequences. We filtered out sequences that contain any unknown nucleotides (other than A, U, C, or G) and sequences whose length exceeds $3,000$nt. This resulted in $45,216$ lncRNA sequences and $79,030$ protein-coding sequences. Next, we applied a clustering algorithm CD-hit~\cite{li2006cd} to lncRNAs and protein-coding RNAs, respectively, to remove redundant sequences. We finally obtained a set of $104,733$ RNA transcripts, consisting of $39,785$ lncRNA sequences and $64,948$ protein-coding sequences.

\subsubsection*{Construction of the lncRNA HTVS pipeline}



For the construction of the lncRNA screening pipeline, we selected four state-of-the-art lncRNA prediction algorithms that have been shown to achieve good prediction performance: CPC$2$~\cite{kang2017CPC2}, CPAT~\cite{wang2013cpat}, PLEK~\cite{li2014plek}, and LncFinder~\cite{han2019lncfinder}.

Table~\ref{tab3} summarizes the performance of the individual algorithm based on the GENCODE dataset, preprocessed as described previously. We assessed the accuracy, sensitivity, and specificity of the respective lncRNA prediction algorithms. For algorithm CPAT, which yields confidence scores between $0$ and $1$ rather than a binary output, we set the decision boundary to $0.5$ for lncRNA classification.
As shown in Table~\ref{tab3}, LncFinder achieved the accuracy, sensitivity, and specificity of $0.8329$, $0.7062$, and $0.9678$, respectively, outperforming all other algorithms in terms of accuracy and sensitivity. However, LncFinder also had the highest computational cost among the compared algorithm, where processing an RNA transcript required $2,495.6231$ milliseconds on average. CPAT was the second-best performer among the four in terms of accuracy and sensitivity. Furthermore, CPAT also achieved the highest specificity. CPC$2$ and PLEK were less accurate compared to LncFinder and CPAT in terms of accuracy, sensitivity, and specificity. Despite their high computational efficiency, both CPC$2$ and CPAT also outperformed PLEK based on overall accuracy.

\begin{figure}[h!]
\centerline{\includegraphics{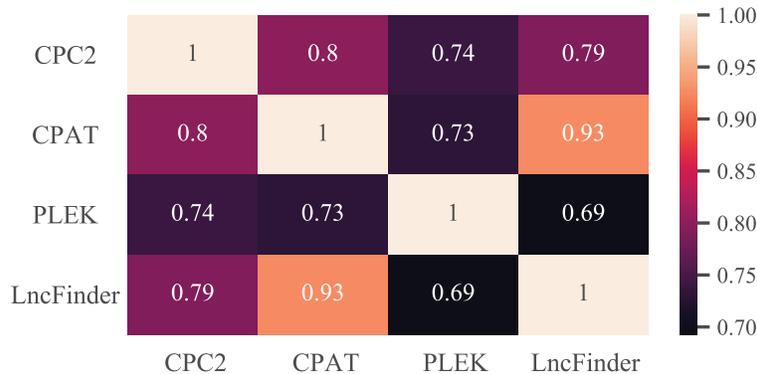}}
\caption{The heat map showing the Pearson's correlation coefficient between different stages. CPAT had the highest correlation to LncFinder. While PLEK was computationally more complex compared to CPAT, it showed a relatively lower correlation to LncFinder.}
\label{fig4}
\end{figure}

\begin{figure}[h!]
\centerline{\includegraphics{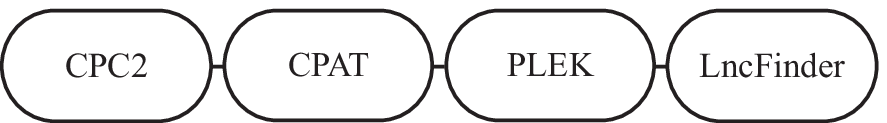}}
\caption{Illustration of a lncRNA HTVS pipeline, where the stages are connected in the order of increasing computational complexity.}
\label{fig3}
\end{figure}

As we previously observed from the performance assessment results based on the synthetic pipeline, the efficacy of the optimized HTVS pipeline is critically dependent on the correlation between the stages constituting the pipeline. The proposed HTVS optimization framework aims to exploit the correlation structure across different screening stages to find the optimal screening policy that strikes the optimal balance between the screening throughput and the computational cost of screening. Here we placed LncFinder--the most accurate and the most computationally costly algorithm among the four--in the final stage. In the first three stages in the HTVS pipeline, we placed CPC$2$, CPAT, and PLEK, in the order of increasing computational complexity. The resulting HTVS pipeline structure is depicted in Fig.~\ref{fig3}. After constructing the screening pipeline, we computed Pearson's correlation coefficient between the predictive output scores obtained from different algorithms.
As shown in Fig.~\ref{fig4}, CPAT showed the highest correlation with LncFinder in the last stage (with a correlation coefficient of $0.93$), the highest among the first three stages in the screening pipeline. 

\begin{table}
\centering
\begin{tabular}{ccccc}
\hline
Algorithm                           & Accuracy & Sensitivity   & Specificity   & Time per RNA (ms) \\ \hline
CPC$2$~\cite{kang2017CPC2}            & 0.7154 & 0.5760        & 0.9493        & 2.5265          \\
CPAT~\cite{wang2013cpat}            & 0.8217 & 0.6861        & 0.9817        & 2.7336          \\
PLEK~\cite{li2014plek}              & 0.7050 & 0.5666        & 0.9478        & 83.1765         \\
LncFinder~\cite{han2019lncfinder}   & 0.8329 & 0.7062        & 0.9678        & 2,495.6231       \\ \hline
\end{tabular}
\caption{Performance of the four individual lncRNA prediction algorithms that constitute the lncRNA HTVS pipeline. The average accuracy, sensitivity, specificity, and processing time (per RNA transcript) are shown.}
\label{tab3}
\end{table}

To apply our proposed HTVS optimization framework, we first estimated the joint probability distribution ${p\left(y_1, y_2, y_3, y_4 \right)}$ of the predictive scores generated by the four different lncRNA prediction algorithms--CPC$2$ ($y_1$), CPAT ($y_2$), PLEK ($y_3$), and LncFinder ($y_4$)--via the EM algorithm~\cite{dempster1977maximum}. For training, $4\%$ of the preprocessed GENCODE data was used. Note that all the computational lncRNA identification algorithms considered in this study output protein-coding probabilities, hence a higher output value corresponds to a higher probability for a given transcript to be protein-coding. Since our goal was to identify the lncRNAs, we multiplied the output scores generated by the algorithms by ${-1}$ such that higher values represent higher chances to be lncRNA transcripts. The screening threshold for the LncFinder in the last stage of the HTVS pipeline was set to ${\lambda_4=0.2}$, which leads to the optimal overall performance of LncFinder with a good balance between sensitivity and specificity.

\subsubsection*{Performance of the optimized lncRNA HTVS pipeline under computational budget constraint}

\begin{figure*}[h!]
\includegraphics[width=\textwidth]{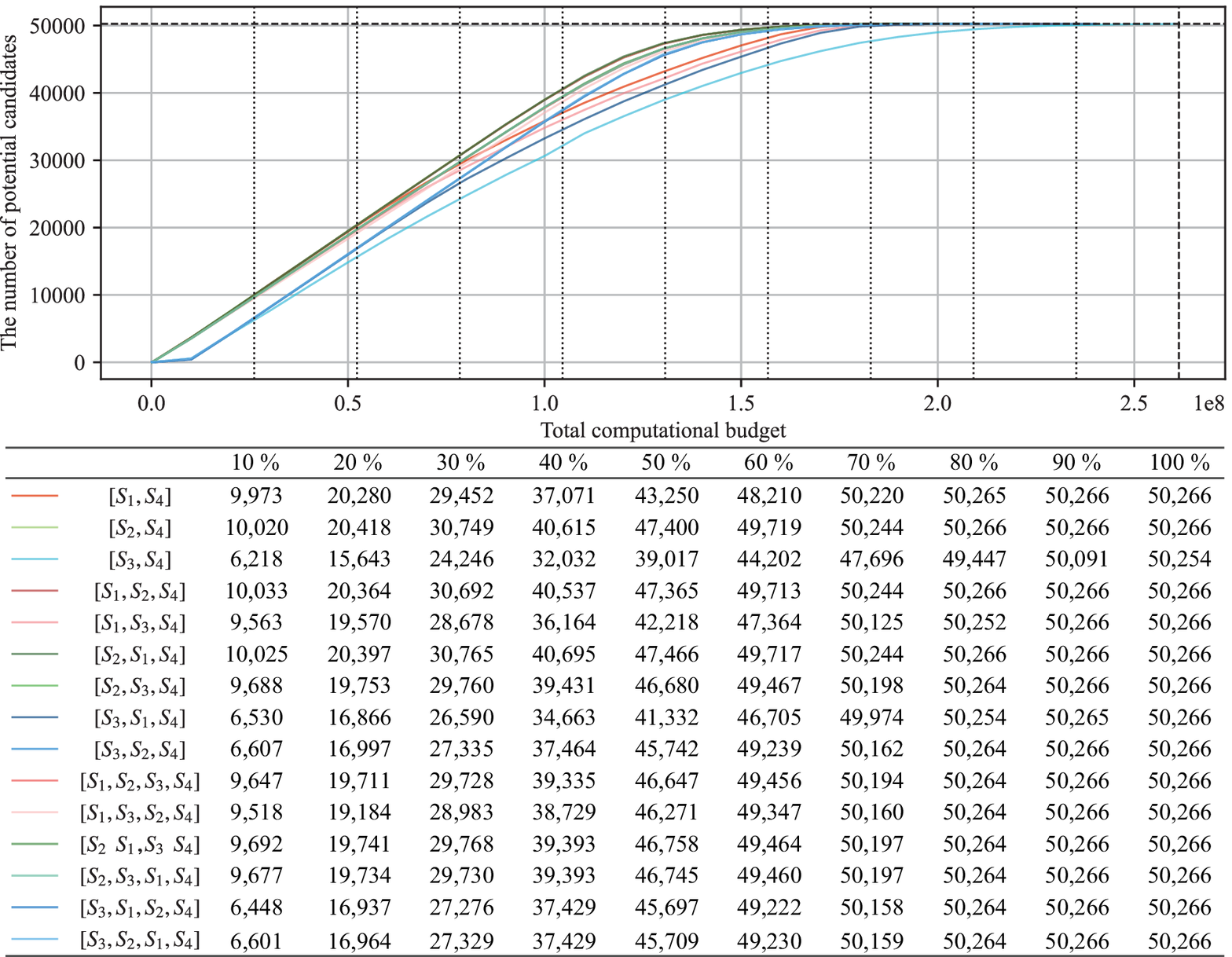}
\caption{Performance evaluation of the optimized lncRNA HTVS pipeline. The number of potential candidates (\textit{i.e.}, lncRNAs) detected by the HTVS pipeline is shown under various computational budget constraints (${x}$-axis). Various different pipeline structures were tested, where the results show that the proposed optimization framework leads to efficient and reliable performance regardless of the structure used.}
\label{fig5}
\end{figure*}

Figure~\ref{fig5} shows the performance of the optimized lncRNA HTVS pipeline for various pipeline structures with the different numbers of stages and ordering. The black horizontal dashed line indicates the total number of potential candidates (\textit{i.e.}, the total number of functional lncRNAs in the test set) and the black vertical dashed line shows the total computational cost (referred to as the ``original cost'' as before) that would be needed for screening all candidates based on the last stage LncFinder alone, without using the HTVS pipeline. Black vertical dotted lines are located at intervals of ${1/10}$ of this original cost. Underneath each dotted line, the number of potential candidates (\textit{i.e.}, true functional lncRNAs) detected by each optimized HTVS pipeline is shown (see the columns in the table aligned with the dotted lines in the plot).

As before, we assumed that the candidates are batch-processed at each stage. As a result, for a given pipeline structure, the computational cost of the first stage determined the minimum computational resources needed to start screening. The correlation between the neighboring stages was closely related to the slope of the corresponding performance curve, which is a phenomenon that we already noticed before based on synthetic pipelines. For example, at ${10\%}$ of the original cost, the pipelines starting with PLEK (\textit{i.e.}, ${S_3}$) showed the worst performance among the tested pipelines in terms of the throughput. Specifically, ${\left[S_3, S_4\right]}$, ${\left[S_3, S_1, S_4\right]}$, ${\left[S_3, S_2, S_4\right]}$, ${\left[S_3, S_1, S_2, S_4\right]}$, and ${\left[S_3, S_2, S_1, S_4\right]}$ detected ${6,218}$; ${6,530}$; ${6,607}$; ${6,448}$; and ${6,601}$ lncRNAs, respectively. On the other hand, pipelines starting with either CPC$2$ or CPAT (\textit{i.e.}, ${S_1}$ or ${S_2}$) detected ${9,518}$ to ${10,033}$ lncRNAs at the same cost. In addition, pipelines ${\left[S_2, S_4\right]}$, ${\left[S_1, S_2, S_4\right]}$, ${\left[S_2, S_1, S_4\right]}$, ${\left[S_2, S_3, S_4\right]}$, ${\left[S_3, S_2, S_4\right]}$, ${\left[S_1,S_2,S_3,S_4\right]}$, ${\left[S_1,S_3, S_2,S_4\right]}$, ${[S_2,}$ ${S_1,S_3,S_4]}$, ${\left[S_2,S_3,S_1,S_4\right]}$, ${\left[S_3,S_1,S_2,S_4\right]}$, and ${\left[S_3,S_2,S_1,S_4\right]}$ including the second stage associated with CPAT that is highly correlated to the last stage LncFinder showed the steepest performance improvement. As a result, all HTVS pipelines that include CPAT were able to identify nearly all true lncRNAs (\textit{i.e.}, $45,697$ to ${47,466}$) at only ${50\%}$ of the original cost, regardless of at which stage CPAT was placed in the pipeline.

While the structure of the HTVS pipeline impacts the overall screening performance, Fig.~\ref{fig5} shows that our proposed optimization framework alleviated the performance dependency on the underlying structure by optimally exploiting the relationships across different stages. For example, even though the optimized pipeline ${\left[S_1, S_2, S_4\right]}$ outperformed the optimized pipeline ${\left[S_1, S_2, S_3, S_4\right]}$, which additionally included PLEK (\textit{i.e.}, ${S_3}$), the performance gap was not very significant. The maximum difference between the two pipeline structures in terms of the detected lncRNAs was ${1,202}$ when the computational budget was set at ${40\%}$ of the original cost. However, when considering that PLEK ($S_3$) was computationally much more expensive compared to CPC$2$ ($S_1$) and CPAT ($S_2$) and also had a lower correlation with LncFinder ($S_4$), the throughput difference of ${1,202}$ was only about ${2.4\%}$ of the total lncRNAs in the test dataset, which is relatively small. Moreover, this throughput difference was drastically reduced as the available computational resources increased. For example, when the computational budget was set at ${70\%}$ of the original cost, the throughput difference between the two pipelines was only ${50}$.

In practice, real-world HTVS pipelines may involve various types of screening stages using multi-fidelity surrogate models. The computational complexity and the fidelity of such surrogate models may differ significantly and the structure of the pipeline may vastly vary depending on the domain experts designing the pipeline. Considering these factors, an important advantage of our proposed HTVS pipeline optimization framework is its capability to consistently attain efficient and accurate screening performance that may weather the effect of potentially suboptimal design choices in constructing real-world HTVS pipelines.

\subsubsection*{Performance of the lncRNA HTVS pipeline jointly optimized for throughput and efficiency}

Next, we evaluated the performance of the lncRNA HTVS pipeline, jointly optimized for both throughput and efficiency based on the proposed framework (with $\alpha = 0.5$). The results for various pipeline configurations are shown in Table~\ref{tab4}.
On average, the optimized HTVS pipeline detected ${48,372}$ lncRNAs out of ${50,266}$ total lncRNAs in the test dataset. The average effective cost was ${3,067}$. Pipeline configurations that include CPAT (${S_2}$) achieved relatively higher computational savings (ranging from ${41.86\%}$ to ${46.25\%}$) compared to those without ${S_2}$ (ranging from ${27.73\%}$ to ${36.52\%}$). As we have previously observed, our proposed optimization framework was effective in maintaining its screening efficiency and accuracy even when the pipeline included a stage (\textit{e.g.}, PLEK) that is less correlated with the last and the highest-fidelity stage (\textit{i.e.}, LncFinder). In fact, the inclusion of a suboptimal stage in the HTVS pipeline does not significantly degrade the average screening performance. This is because the proposed optimization framework enables one to select the optimal threshold values that can sensibly combine the benefits of the most efficient stages (such as CPC$2$ and CPAT in this case) as well as the most accurate stage (LncFinder), thereby maximizing the expected ROCI.
Similar observation can be made regarding the ordering of the multiple screening stages, as Table~\ref{tab4} shows that the average performance does not significantly depend on the actual ordering of the stages when the screening threshold values are optimized via our proposed framework.
For example, when using all four stages in the HTVS pipeline (${N = 4}$), the optimized pipeline detected ${48,402}$ lncRNAs on average and with consistent computational savings ranging between ${42.27\%}$ and ${44.00\%}$. 
We also evaluated the accuracy of the potential candidates screened by the optimized HTVS pipeline based on four performance metrics: \textit{accuracy (ACC)}, \textit{sensitivity (SN)}, \textit{specificity (SP)}, and \textit{F1} score. Interestingly, all configurations except for ${\left[S_1, S_4\right]}$ outperformed LncFinder in terms of ACC. In terms of SN, the optimized pipeline achieved an average sensitivity of $0.9177$. All pipeline configurations resulted in higher specificity compared to LncFinder. Besides, pipeline configurations that include ${S_2}$ consistently outperformed LncFinder in terms of the F1 score.

\begin{table*}
\scriptsize
\centering
\begin{tabular}{ccccccccc}
\hline
Configuration     & \begin{tabular}[c]{@{}c@{}}Potential\\ candidates\end{tabular} & Total cost (ms) & \begin{tabular}[c]{@{}c@{}}Effective\\cost\end{tabular} & \begin{tabular}[c]{@{}c@{}}Computational\\savings\end{tabular} & Accuracy & Sensitivity & Specificity & F1 \\ \hline
{${\left[S_4\right]}$}               & $50,266$                       & ${261,374,090}$ & ${5,200}$ & ${0\%}$ & ${0.8440}$ & ${0.9264}$ & ${0.7936}$ & ${0.8186}$ \\ \hline
{${\left[S_1, S_4\right]}$}          & {${48,875}$}    & ${161,357,081}$ & ${3,301}$ & ${36.52\%}$ & ${0.8429}$ & ${0.9075}$ & ${0.8034}$ & ${0.8144}$ \\
{${\left[S_2, S_4\right]}$}          & {${47,950}$}    & ${134,366,143}$ & ${2,802}$ & ${46.12\%}$ & ${0.8624}$ & ${0.9215}$ & ${0.8262}$ & ${0.8357}$ \\
{${\left[S_3, S_4\right]}$}          & {${47,083}$}    & ${176,963,736}$ & ${3,758}$ & ${27.73\%}$ & ${0.8450}$ & ${0.8876}$ & ${0.8188}$ & ${0.8131}$ \\
{${\left[S_1, S_2, S_4\right]}$}     & {${48,210}$}    & ${134,748,992}$ & ${2,795}$ & ${46.25\%}$ & ${0.8600}$ & ${0.9216}$ & ${0.8222}$ & ${0.8333}$ \\
{${\left[S_1, S_3, S_4\right]}$}     & {${49,100}$}    & ${168,490,516}$ & ${3,432}$ & ${34.00\%}$ & ${0.8442}$ & ${0.9120}$ & ${0.8026}$ & ${0.8164}$ \\
{${\left[S_2, S_1, S_4\right]}$}     & {${48,214}$}    & ${134,812,024}$ & ${2,796}$ & ${46.23\%}$ & ${0.8600}$ & ${0.9216}$ & ${0.8222}$ & ${0.8334}$ \\
{${\left[S_2, S_3, S_4\right]}$}     & {${48,295}$}    & ${141,710,246}$ & ${2,934}$ & ${43.58\%}$ & ${0.8602}$ & ${0.9230}$ & ${0.8218}$ & ${0.8338}$ \\
{${\left[S_3, S_1, S_4\right]}$}     & {${49,119}$}    & ${171,803,403}$ & ${3,498}$ & ${32.73\%}$ & ${0.8444}$ & ${0.9124}$ & ${0.8026}$ & ${0.8166}$ \\
{${\left[S_3, S_2, S_4\right]}$}     & {${48,326}$}    & ${146,100,080}$ & ${3,023}$ & ${41.86\%}$ & ${0.8600}$ & ${0.9231}$ & ${0.8214}$ & ${0.8336}$ \\
{${\left[S_1, S_2, S_3, S_4\right]}$}    & {${48,402}$}& ${140,954,256}$ & ${2,912}$ & ${44.00\%}$ & ${0.8591}$ & ${0.9228}$ & ${0.8200}$ & ${0.8326}$ \\
{${\left[S_1, S_3, S_2, S_4\right]}$}    & {${48,332}$}& ${141,229,518}$ & ${2,922}$ & ${43.81\%}$ & ${0.8587}$ & ${0.9215}$ & ${0.8203}$ & ${0.8321}$ \\
{${\left[S_2, S_1, S_3, S_4\right]}$}    & {${48,409}$}& ${141,022,859}$ & ${2,913}$ & ${43.98\%}$ & ${0.8591}$ & ${0.9229}$ & ${0.8200}$ & ${0.8326}$ \\
{${\left[S_2, S_3, S_1, S_4\right]}$}    & {${48,414}$}& ${141,225,328}$ & ${2,917}$ & ${43.90\%}$ & ${0.8591}$ & ${0.9230}$ & ${0.8200}$ & ${0.8327}$ \\
{${\left[S_3, S_1, S_2, S_4\right]}$}    & {${48,424}$}& ${145,321,388}$ & ${3,001}$ & ${42.29\%}$ & ${0.8589}$ & ${0.9228}$ & ${0.8197}$ & ${0.8324}$ \\
{${\left[S_3, S_2, S_1, S_4\right]}$}    & {${48,429}$}& ${145,388,626}$ & ${3,002}$ & ${42.27\%}$ & ${0.8589}$ & ${0.9229}$ & ${0.8197}$ & ${0.8325}$ \\\hline
\end{tabular}
\caption{Performance evaluation of the lncRNA HTVS pipeline jointly optimized for throughput and efficiency. Results are shown for various pipeline configurations, where the optimized screening policy was used (with ${\alpha = 0.5}$).}
\label{tab4}
\end{table*}

Finally, we compared the performance of the optimized pipeline to that of the baseline approach that selects the top $R_s \%$ of the incoming candidates for the next stage, where $R_s \%$ is a parameter to be determined by a domain expert. For this comparison, we considered the four-stage pipeline ${\left[S_1, S_2, S_3, S_4\right]}$. The optimal screening policy was found based on our proposed framework using three different values of $\alpha \in \left\{0.25, 0.50, 0.75\right\}$. The baseline screening approach was evaluated based on four different levels of $R_s \in \left\{25\%, 50\%, 75\%\right\}$. The performance assessment results are summarized in Table~\ref{tab5}.
As shown in Table~\ref{tab5}, the baseline approach detected fewer lncRNAs for all values of $R_s$ compared to the optimized pipeline. Specifically, the jointly optimized pipeline detected ${48,965}$; ${48,402}$; and ${47,106}$ lncRNAs at a cost of ${148,155,016}$ ($\alpha=0.75$); ${140,954,256}$ ($\alpha=0.50$); and ${131,830,857}$ ($\alpha=0.25$), respectively. On the other hand, the baseline approach with $R_s = 75\%$ detected only ${39,079}$ lncRNAs at a total cost of ${115,643,459}$. For $R_s = 50\%$ and $R_s = 25\%$, the baseline scheme detected only ${12,653}$ and ${1,402}$ lncRNAs, respectively. In terms of the four quality metrics (ACC, SN, SP, and F1), the optimized pipeline outperformed the baseline scheme in terms of ACC, SN, and F1. The optimized pipeline resulted in lower SP compared to the baseline. However, it should be noted that the potential candidates detected by the optimized HTVS pipeline are remarkably higher compared to the baseline approach. This is clearly reflected in the much lower SN of the baseline approach, as shown in Table~\ref{tab5}. As a result, the baseline approach tends to achieve significantly lower ACC and F1 compared to the optimal screening scheme.

\begin{table*}
\scriptsize
\centering
\begin{tabular}{ccccccccc}
\hline
Approach     & \begin{tabular}[c]{@{}c@{}}Potential\\ candidates\end{tabular} & Total cost (ms) & \begin{tabular}[c]{@{}c@{}}Effective\\cost\end{tabular} & \begin{tabular}[c]{@{}c@{}}Computational\\savings\end{tabular} & Accuracy & Sensitivity & Specificity & F1 \\ \hline
Proposed (${\alpha = 0.75}$) & {${48,965}$}& ${148,155,016}$ & ${3,026}$ & ${41.81\%}$ & ${0.8553}$ & ${0.9249}$ & ${0.8126}$ & ${0.8292}$ \\
Proposed (${\alpha = 0.5}$) & {${48,402}$}& ${140,954,256}$ & ${2,912}$ & ${44.00\%}$ & ${0.8591}$ & ${0.9228}$ & ${0.8200}$ & ${0.8326}$ \\
Proposed (${\alpha = 0.25}$) & {${47,106}$}& ${131,830,857}$ & ${2,799}$ & ${46.17\%}$ & ${0.8650}$ & ${0.9143}$ & ${0.8348}$ & ${0.8373}$ \\ \hline
Baseline ($R_s = 0.75$)     & {${39,079}$}& ${115,643,459}$ & ${2,959}$ & ${43.10\%}$ & ${0.8366}$ & ${0.7761}$ & ${0.8737}$ & ${0.7801}$ \\
Baseline ($R_s = 0.5$)      & {${12,653}$}& ${35,255,772}$  & ${2,786}$ & ${46.42\%}$ & ${0.7170}$ & ${0.2866}$ & ${0.9807}$ & ${0.4348}$ \\ 
Baseline ($R_s = 0.25$)     & {${1,402}$}&  ${4,963,415}$   & ${3,540}$ & ${31.92\%}$ & ${0.6318}$ & ${0.0330}$ & ${0.9986}$ & ${0.0638}$ \\ \hline
\end{tabular}
\caption{Performance of the four-stage lncRNA HTVS pipeline {${\left[S_1, S_2, S_3, S_4\right]}$}. The overall performance of the HTVS pipeline jointly optimized for throughput and efficiency is compared to that of the baseline screening approach.}
\label{tab5}
\end{table*}

\section*{DISCUSSION}
\label{sec:discussion}

In this work, we proposed a general mathematical framework for identifying the optimal screening policy that can maximize the return-on-computational-investment (ROCI) of an HTVS pipeline. The need for screening a large set of molecules to detect potential candidates that possess the desired properties frequently arise in various science and engineering domains, although the design and operation of such screening pipelines strongly depend on expert intuitions and \textit{ad hoc} approaches. We aimed to rectify this problem by taking a principled approach to high-throughput virtual screening (HTVS), thereby maximizing the screening performance of a given HTVS pipeline, reducing the performance dependence on the pipeline configuration, and enabling quantitative comparison between different HTVS pipelines based on their optimal achievable performance.

We considered two scenarios for HTVS performance optimization in this study: first, maximizing the detection of potential candidate molecules that possess the desired property under a constrained computational budget; second, jointly optimizing the throughput and the computational efficiency of the HTVS pipeline when there is no fixed computational budget for the screening operation. 
For both scenarios, we have thoroughly tested the performance of our proposed HTVS optimization framework. Comprehensive performance assessment based on synthetic HTVS pipelines as well as real lncRNA screening pipelines both showed clear advantages of the proposed framework. 
Not only does the HTVS optimization framework remove the guesswork in the operation of HTVS pipelines to maximize the throughput, enhance the screening accuracy, and minimize the computational cost, it leads to reliable and consistent screening performance across a wide variety of HTVS pipeline structures. This is a significant benefit of the proposed framework that is of practical importance since it makes the overall screening performance robust to variations and potentially suboptimal design choices in constructing real-world HTVS pipelines. As there can be infinite different ways of building an HTVS pipeline in real scientific and engineering applications, it is important to note that our proposed optimization framework can guarantee near-optimal screening performance for any reasonable design choice regarding the HTVS pipeline configuration. 

While the HTVS optimization framework presented in this paper is fairly general and may be applied to various molecular screening problems, there are interesting open research problems for extending its capabilities even further. For example, how can we design an optimal screening policy when the structure of the HTVS pipeline is not linear (\textit{i.e.}, sequential)? For example, it would be interesting to investigate how one may design optimal screening policies for pipelines whose structure is a directed acyclic graph (DAG). Another interesting problem is how to construct an ideal screening pipeline if only a high-fidelity model is available and its lower-fidelity surrogates need to be designed/learned from this high-fidelity model from scratch. Finally, we can envision designing a dynamic screening policy--\textit{e.g.}, based on reinforcement learning--where there is no predetermined ``pipeline structure'' but a candidate molecule is dynamically steered to the next stage based on the screening outcome of the previous stage. In this case, instead of first learning the relations between different stages and then deriving the optimal screening policy--as we did in this current study--we may learn, apply, and enhance the screening policy in a dynamic and adaptive manner, which may be better suited for certain applications.

\bibliographystyle{ieeetr}
\bibliography{main.bib}

\begin{thebibliography}{10}

\bibitem{saadi2020impeccable}
A.~A. Saadi, D.~Alfe, Y.~Babuji, A.~Bhati, B.~Blaiszik, A.~Brace, T.~Brettin,
  K.~Chard, R.~Chard, A.~Clyde, {\em et~al.}, ``Impeccable: integrated modeling
  pipeline for covid cure by assessing better leads,'' in {\em 50th
  International Conference on Parallel Processing}, pp.~1--12, 2021.

\bibitem{roy2021global}
B.~Roy, J.~Dhillon, N.~Habib, and B.~Pugazhandhi, ``Global variants of
  covid-19: Current understanding,'' {\em Journal of Biomedical Sciences},
  vol.~8, no.~1, pp.~8--11, 2021.

\bibitem{sterling2015zinc}
T.~Sterling and J.~J. Irwin, ``Zinc 15--ligand discovery for everyone,'' {\em
  Journal of chemical information and modeling}, vol.~55, no.~11,
  pp.~2324--2337, 2015.

\bibitem{rieber2009rnaither}
N.~Rieber, B.~Knapp, R.~Eils, and L.~Kaderali, ``Rnaither, an automated
  pipeline for the statistical analysis of high-throughput rnai screens,'' {\em
  Bioinformatics}, vol.~25, no.~5, pp.~678--679, 2009.

\bibitem{studer2010engineering}
M.~H. Studer, J.~D. DeMartini, S.~Brethauer, H.~L. McKenzie, and C.~E. Wyman,
  ``Engineering of a high-throughput screening system to identify cellulosic
  biomass, pretreatments, and enzyme formulations that enhance sugar release,''
  {\em Biotechnology and Bioengineering}, vol.~105, no.~2, pp.~231--238, 2010.

\bibitem{hartmann2011htpheno}
A.~Hartmann, T.~Czauderna, R.~Hoffmann, N.~Stein, and F.~Schreiber, ``Htpheno:
  an image analysis pipeline for high-throughput plant phenotyping,'' {\em BMC
  bioinformatics}, vol.~12, no.~1, pp.~1--9, 2011.

\bibitem{sikorski2018high}
K.~Sikorski, A.~Mehta, M.~Inngjerdingen, F.~Thakor, S.~Kling, T.~Kalina, T.~A.
  Nyman, M.~E. Stensland, W.~Zhou, G.~A. de~Souza, {\em et~al.}, ``A
  high-throughput pipeline for validation of antibodies,'' {\em Nature
  methods}, vol.~15, no.~11, pp.~909--912, 2018.

\bibitem{clyde2021high}
A.~Clyde, S.~Galanie, D.~W. Kneller, H.~Ma, Y.~Babuji, B.~Blaiszik, A.~Brace,
  T.~Brettin, K.~Chard, R.~Chard, {\em et~al.}, ``High-throughput virtual
  screening and validation of a sars-cov-2 main protease noncovalent
  inhibitor,'' {\em Journal of chemical information and modeling}, vol.~62,
  no.~1, pp.~116--128, 2021.

\bibitem{martin2014silico}
R.~L. Martin, C.~M. Simon, B.~Smit, and M.~Haranczyk, ``In silico design of
  porous polymer networks: high-throughput screening for methane storage
  materials,'' {\em Journal of the American Chemical Society}, vol.~136,
  no.~13, pp.~5006--5022, 2014.

\bibitem{cheng2015accelerating}
L.~Cheng, R.~S. Assary, X.~Qu, A.~Jain, S.~P. Ong, N.~N. Rajput, K.~Persson,
  and L.~A. Curtiss, ``Accelerating electrolyte discovery for energy storage
  with high-throughput screening,'' {\em The journal of physical chemistry
  letters}, vol.~6, no.~2, pp.~283--291, 2015.

\bibitem{chen2017developing}
J.~J.~F. Chen and D.~P. Visco~Jr, ``Developing an in silico pipeline for faster
  drug candidate discovery: Virtual high throughput screening with the
  signature molecular descriptor using support vector machine models,'' {\em
  Chemical Engineering Science}, vol.~159, pp.~31--42, 2017.

\bibitem{filer2017tcpl}
D.~L. Filer, P.~Kothiya, R.~W. Setzer, R.~S. Judson, and M.~T. Martin, ``tcpl:
  the toxcast pipeline for high-throughput screening data,'' {\em
  Bioinformatics}, vol.~33, no.~4, pp.~618--620, 2017.

\bibitem{rebbeck2020ryr1}
R.~T. Rebbeck, D.~P. Singh, K.~A. Janicek, D.~M. Bers, D.~D. Thomas, B.~S.
  Launikonis, and R.~L. Cornea, ``Ryr1-targeted drug discovery pipeline
  integrating fret-based high-throughput screening and human myofiber dynamic
  ca 2+ assays,'' {\em Scientific reports}, vol.~10, no.~1, pp.~1--13, 2020.

\bibitem{tran2020smf}
A.~Tran, T.~Wildey, and S.~McCann, ``smf-bo-2cogp: A sequential multi-fidelity
  constrained bayesian optimization framework for design applications,'' {\em
  Journal of Computing and Information Science in Engineering}, vol.~20, no.~3,
  2020.

\bibitem{yan2017solar}
Q.~Yan, J.~Yu, S.~K. Suram, L.~Zhou, A.~Shinde, P.~F. Newhouse, W.~Chen, G.~Li,
  K.~A. Persson, J.~M. Gregoire, {\em et~al.}, ``Solar fuels photoanode
  materials discovery by integrating high-throughput theory and experiment,''
  {\em Proceedings of the National Academy of Sciences}, vol.~114, no.~12,
  pp.~3040--3043, 2017.

\bibitem{zhang2020first}
B.~Zhang, X.~Zhang, J.~Yu, Y.~Wang, K.~Wu, and M.-H. Lee, ``First-principles
  high-throughput screening pipeline for nonlinear optical materials:
  Application to borates,'' {\em Chemistry of Materials}, vol.~32, no.~15,
  pp.~6772--6779, 2020.

\bibitem{chen1995new}
C.~Chen, Y.~Wang, Y.~Xia, B.~Wu, D.~Tang, K.~Wu, Z.~Wenrong, L.~Yu, and L.~Mei,
  ``New development of nonlinear optical crystals for the ultraviolet region
  with molecular engineering approach,'' {\em Journal of applied physics},
  vol.~77, no.~6, pp.~2268--2272, 1995.

\bibitem{shi2017finding}
G.~Shi, Y.~Wang, F.~Zhang, B.~Zhang, Z.~Yang, X.~Hou, S.~Pan, and K.~R.
  Poeppelmeier, ``Finding the next deep-ultraviolet nonlinear optical material:
  Nh4b4o6f,'' {\em Journal of the American Chemical Society}, vol.~139, no.~31,
  pp.~10645--10648, 2017.

\bibitem{zhang2017fluorooxoborates}
B.~Zhang, G.~Shi, Z.~Yang, F.~Zhang, and S.~Pan, ``Fluorooxoborates:
  beryllium-free deep-ultraviolet nonlinear optical materials without layered
  growth,'' {\em Angewandte Chemie International Edition}, vol.~56, no.~14,
  pp.~3916--3919, 2017.

\bibitem{luo2018m2b10o14f6}
M.~Luo, F.~Liang, Y.~Song, D.~Zhao, F.~Xu, N.~Ye, and Z.~Lin, ``M2b10o14f6 (m=
  ca, sr): Two noncentrosymmetric alkaline earth fluorooxoborates as promising
  next-generation deep-ultraviolet nonlinear optical materials,'' {\em Journal
  of the American Chemical Society}, vol.~140, no.~11, pp.~3884--3887, 2018.

\bibitem{mutailipu2018srb5o7f3}
M.~Mutailipu, M.~Zhang, B.~Zhang, L.~Wang, Z.~Yang, X.~Zhou, and S.~Pan,
  ``Srb5o7f3 functionalized with [b5o9f3] 6- chromophores: Accelerating the
  rational design of deep-ultraviolet nonlinear optical materials,'' {\em
  Angewandte Chemie}, vol.~130, no.~21, pp.~6203--6207, 2018.

\bibitem{wang2018cation}
Y.~Wang, B.~Zhang, Z.~Yang, and S.~Pan, ``Cation-tuned synthesis of
  fluorooxoborates: Towards optimal deep-ultraviolet nonlinear optical
  materials,'' {\em Angewandte Chemie}, vol.~130, no.~8, pp.~2172--2176, 2018.

\bibitem{zhang2018cab5o7f3}
Z.~Zhang, Y.~Wang, B.~Zhang, Z.~Yang, and S.~Pan, ``Cab5o7f3: A beryllium-free
  alkaline-earth fluorooxoborate exhibiting excellent nonlinear optical
  performances,'' {\em Inorganic chemistry}, vol.~57, no.~9, pp.~4820--4823,
  2018.

\bibitem{wilmer2012large}
C.~E. Wilmer, M.~Leaf, C.~Y. Lee, O.~K. Farha, B.~G. Hauser, J.~T. Hupp, and
  R.~Q. Snurr, ``Large-scale screening of hypothetical metal--organic
  frameworks,'' {\em Nature chemistry}, vol.~4, no.~2, pp.~83--89, 2012.

\bibitem{gupta2020computational}
S.~Gupta, D.~Parihar, M.~Shah, S.~Yadav, H.~Managori, S.~Bhowmick, P.~C. Patil,
  S.~A. Alissa, S.~M. Wabaidur, and M.~A. Islam, ``Computational screening of
  promising beta-secretase 1 inhibitors through multi-step molecular docking
  and molecular dynamics simulations-pharmacoinformatics approach,'' {\em
  Journal of Molecular Structure}, vol.~1205, p.~127660, 2020.

\bibitem{kim2021universal}
D.~Y. Kim, M.~Ha, and K.~S. Kim, ``A universal screening strategy for the
  accelerated design of superior oxygen evolution/reduction electrocatalysts,''
  {\em Journal of Materials Chemistry A}, vol.~9, no.~6, pp.~3511--3519, 2021.

\bibitem{liu2015high}
T.~Liu, K.~C. Kim, R.~Kavian, S.~S. Jang, and S.~W. Lee, ``High-density
  lithium-ion energy storage utilizing the surface redox reactions in folded
  graphene films,'' {\em Chemistry of Materials}, vol.~27, no.~9,
  pp.~3291--3298, 2015.

\bibitem{kim2016first}
K.~C. Kim, T.~Liu, S.~W. Lee, and S.~S. Jang, ``First-principles density
  functional theory modeling of li binding: thermodynamics and redox properties
  of quinone derivatives for lithium-ion batteries,'' {\em Journal of the
  American Chemical Society}, vol.~138, no.~7, pp.~2374--2382, 2016.

\bibitem{kim2016thermodynamic}
S.~Kim, K.~C. Kim, S.~W. Lee, and S.~S. Jang, ``Thermodynamic and redox
  properties of graphene oxides for lithium-ion battery applications: a first
  principles density functional theory modeling approach,'' {\em Physical
  Chemistry Chemical Physics}, vol.~18, no.~30, pp.~20600--20606, 2016.

\bibitem{liu2017self}
T.~Liu, K.~C. Kim, B.~Lee, Z.~Chen, S.~Noda, S.~S. Jang, and S.~W. Lee,
  ``Self-polymerized dopamine as an organic cathode for li-and na-ion
  batteries,'' {\em Energy \& Environmental Science}, vol.~10, no.~1,
  pp.~205--215, 2017.

\bibitem{park2017systematic}
J.~H. Park, T.~Liu, K.~C. Kim, S.~W. Lee, and S.~S. Jang, ``Systematic
  molecular design of ketone derivatives of aromatic molecules for lithium-ion
  batteries: First-principles dft modeling,'' {\em ChemSusChem}, vol.~10,
  no.~7, pp.~1584--1591, 2017.

\bibitem{kang2018density}
J.~Kang, K.~C. Kim, and S.~S. Jang, ``Density functional theory
  modeling-assisted investigation of thermodynamics and redox properties of
  boron-doped corannulenes for cathodes in lithium-ion batteries,'' {\em The
  Journal of Physical Chemistry C}, vol.~122, no.~20, pp.~10675--10681, 2018.

\bibitem{sood2018electrochemical}
P.~Sood, K.~C. Kim, and S.~S. Jang, ``Electrochemical and electronic properties
  of nitrogen doped fullerene and its derivatives for lithium-ion battery
  applications,'' {\em Journal of energy chemistry}, vol.~27, no.~2,
  pp.~528--534, 2018.

\bibitem{sood2018electrochemical2}
P.~Sood, K.~C. Kim, and S.~S. Jang, ``Electrochemical properties of boron-doped
  fullerene derivatives for lithium-ion battery applications,'' {\em
  ChemPhysChem}, vol.~19, no.~6, pp.~753--758, 2018.

\bibitem{zhu2018boron}
Y.~Zhu, K.~C. Kim, and S.~S. Jang, ``Boron-doped coronenes with high redox
  potential for organic positive electrodes in lithium-ion batteries: a
  first-principles density functional theory modeling study,'' {\em Journal of
  Materials Chemistry A}, vol.~6, no.~21, pp.~10111--10120, 2018.

\bibitem{dempster1977maximum}
A.~P. Dempster, N.~M. Laird, and D.~B. Rubin, ``Maximum likelihood from
  incomplete data via the em algorithm,'' {\em Journal of the Royal Statistical
  Society: Series B (Methodological)}, vol.~39, no.~1, pp.~1--22, 1977.

\bibitem{ng2013long}
S.-Y. Ng, L.~Lin, B.~S. Soh, and L.~W. Stanton, ``Long noncoding rnas in
  development and disease of the central nervous system,'' {\em Trends in
  Genetics}, vol.~29, no.~8, pp.~461--468, 2013.

\bibitem{tan2013non}
L.~Tan, J.-T. Yu, N.~Hu, and L.~Tan, ``Non-coding rnas in alzheimer's
  disease,'' {\em Molecular neurobiology}, vol.~47, no.~1, pp.~382--393, 2013.

\bibitem{luo2016long}
Q.~Luo and Y.~Chen, ``Long noncoding rnas and alzheimer’s disease,'' {\em
  Clinical interventions in aging}, vol.~11, p.~867, 2016.

\bibitem{congrains2012genetic}
A.~Congrains, K.~Kamide, R.~Oguro, O.~Yasuda, K.~Miyata, E.~Yamamoto, T.~Kawai,
  H.~Kusunoki, H.~Yamamoto, Y.~Takeya, {\em et~al.}, ``Genetic variants at the
  9p21 locus contribute to atherosclerosis through modulation of anril and
  cdkn2a/b,'' {\em Atherosclerosis}, vol.~220, no.~2, pp.~449--455, 2012.

\bibitem{xue2016g}
Z.~Xue, S.~Hennelly, B.~Doyle, A.~A. Gulati, I.~V. Novikova, K.~Y. Sanbonmatsu,
  and L.~A. Boyer, ``A g-rich motif in the lncrna braveheart interacts with a
  zinc-finger transcription factor to specify the cardiovascular lineage,''
  {\em Molecular cell}, vol.~64, no.~1, pp.~37--50, 2016.

\bibitem{yang2014lncrna}
G.~Yang, X.~Lu, and L.~Yuan, ``Lncrna: a link between rna and cancer,'' {\em
  Biochimica et Biophysica Acta (BBA)-Gene Regulatory Mechanisms}, vol.~1839,
  no.~11, pp.~1097--1109, 2014.

\bibitem{shi2015critical}
X.~Shi, M.~Sun, H.~Liu, Y.~Yao, R.~Kong, F.~Chen, and Y.~Song, ``A critical
  role for the long non-coding rna gas5 in proliferation and apoptosis in
  non-small-cell lung cancer,'' {\em Molecular carcinogenesis}, vol.~54,
  no.~S1, pp.~E1--E12, 2015.

\bibitem{peng2017lncrna}
W.-X. Peng, P.~Koirala, and Y.-Y. Mo, ``Lncrna-mediated regulation of cell
  signaling in cancer,'' {\em Oncogene}, vol.~36, no.~41, pp.~5661--5667, 2017.

\bibitem{carlevaro2020cancer}
J.~Carlevaro-Fita, A.~Lanz{\'o}s, L.~Feuerbach, C.~Hong, D.~Mas-Ponte, J.~S.
  Pedersen, and R.~Johnson, ``Cancer lncrna census reveals evidence for deep
  functional conservation of long noncoding rnas in tumorigenesis,'' {\em
  Communications biology}, vol.~3, no.~1, pp.~1--16, 2020.

\bibitem{wang2013cpat}
L.~Wang, H.~J. Park, S.~Dasari, S.~Wang, J.-P. Kocher, and W.~Li, ``Cpat:
  Coding-potential assessment tool using an alignment-free logistic regression
  model,'' {\em Nucleic acids research}, vol.~41, no.~6, pp.~e74--e74, 2013.

\bibitem{li2014plek}
A.~Li, J.~Zhang, and Z.~Zhou, ``Plek: a tool for predicting long non-coding
  rnas and messenger rnas based on an improved k-mer scheme,'' {\em BMC
  bioinformatics}, vol.~15, no.~1, pp.~1--10, 2014.

\bibitem{kang2017CPC2}
Y.-J. Kang, D.-C. Yang, L.~Kong, M.~Hou, Y.-Q. Meng, L.~Wei, and G.~Gao,
  ``Cpc2: a fast and accurate coding potential calculator based on sequence
  intrinsic features,'' {\em Nucleic acids research}, vol.~45, no.~W1,
  pp.~W12--W16, 2017.

\bibitem{han2019lncfinder}
S.~Han, Y.~Liang, Q.~Ma, Y.~Xu, Y.~Zhang, W.~Du, C.~Wang, and Y.~Li,
  ``Lncfinder: an integrated platform for long non-coding rna identification
  utilizing sequence intrinsic composition, structural information and
  physicochemical property,'' {\em Briefings in bioinformatics}, vol.~20,
  no.~6, pp.~2009--2027, 2019.

\bibitem{frankish2021gencode}
A.~Frankish, M.~Diekhans, I.~Jungreis, J.~Lagarde, J.~E. Loveland, J.~M. Mudge,
  C.~Sisu, J.~C. Wright, J.~Armstrong, I.~Barnes, {\em et~al.}, ``Gencode
  2021,'' {\em Nucleic acids research}, vol.~49, no.~D1, pp.~D916--D923, 2021.

\bibitem{li2006cd}
W.~Li and A.~Godzik, ``Cd-hit: a fast program for clustering and comparing
  large sets of protein or nucleotide sequences,'' {\em Bioinformatics},
  vol.~22, no.~13, pp.~1658--1659, 2006.

\end{thebibliography}
\end{document}